\documentclass[final]{siamltex}
\usepackage{amssymb,amstext,mathrsfs}
\usepackage[leqno]{amsmath}
\usepackage{cite}
\usepackage{float}
\usepackage{graphics}
\usepackage{graphicx}
\usepackage{multirow}
\usepackage[colorlinks,linkcolor=red,anchorcolor=black,citecolor=blue]{hyperref}
\usepackage[figure,table]{hypcap}

\newtheorem{scheme}{Scheme}
\newtheorem{remark}[theorem]{Remark}
\newtheorem{example}{{\bf Example}}
\newtheorem{Assum}{Assumption}
\newtheorem{thm}{Theorem}

\def\F{\mathscr{F}}

\def\R{\mathbb{R}}

\def\11{\mathbf{1}}
\def\aa{\mathbf{a}}

\def\jj{\mathbf{j}}
\def\xx{\mathbf{x}}

\def\di{\mathrm{d}}
\newcommand{\abs}[1]{\left\vert#1\right\vert}

\newcommand{\Ec}[3]{\mathbb{E}_{#1}^{#2}\left[#3\right]}

\newcommand{\CE}[2]{\mathbb{E}\left[\left.#1\right\vert #2\right]}

\title{A New Kind of High-Order Multistep Schemes
for Coupled Forward Backward Stochastic Differential Equations
\thanks{This work is partially supported  by  the National
Natural Science Foundations of China under grant numbers 91130003, 11201461and 11171189, and
Natural Science Foundation of Shandong Province under grant number ZR2011AZ002.}
}

\author{Weidong Zhao\thanks{School of Mathematics \& Finance Institute, Shandong University, Jinan 250100, China.
Email: wdzhao@sdu.edu.cn.}
        \and Yu Fu\thanks{School of Mathematics \& Finance Institute, Shandong University, Jinan 250100, China.
Email: nielf0614@126.com.}
\and Tao Zhou\thanks{LSEC, Institute
of Computational Mathematics, Academy of Mathematics and Systems
Science, Chinese Academy of Sciences, Beijing 100190,
China. Email: tzhou@lsec.cc.ac.cn.}
}
\begin{document}

\maketitle

\begin{abstract}
In this work, we are concerned with the high-order numerical methods for coupled forward-backward stochastic differential
equations (FBSDEs). Based on the FBSDEs theory,
we derive two reference ordinary differential equations (ODEs) from the backward SDE, which contain the conditional expectations and their derivatives. Then, our high-order multistep schemes are obtained
by carefully approximating the conditional
expectations and the derivatives, in the reference ODEs. Motivated by the local property of the generator of diffusion processes, the Euler method is used to solve the \textit{forward} SDE, however, it is noticed that the numerical solution of the \textit{backward} SDE is still of high-order accuracy. Such results are obviously promising: on one hand, the use of Euler method (for the forward SDE) can dramatically simplify the entire computational scheme, and on the other hand, one might be only interested in the solution of the backward SDE in many real applications such as option pricing. Several numerical experiments
are presented to demonstrate the effectiveness of the numerical method.

\end{abstract}



\pagestyle{myheadings}
\thispagestyle{plain}
\markboth{W. Zhao, Y. Fu and T. Zhou}
{High-Order Multistep Schemes
for Coupled FBSDEs}

\section{Introduction}

We are interested in the numerical solution of the following coupled forward-backward stochastic differential equations~(FBSDEs)
on $(\Omega,\mathcal{F}, \mathbb F, P)$.
\begin{equation}\label{FBSDEs}
\left\{
\begin{aligned}
X_t &= X_0 + \int_0^t b(s,X_s,Y_s,Z_s) \di t + \int_0^t \sigma(s,X_s,Y_s,Z_s)\di W_s,\,\, &{\rm SDE}\\
Y_t &= \xi + \int_t^T f(s,X_s,Y_s,Z_s)\di s - \int_t^T Z_s\di W_s, \,\, &{\rm BSDE}
\end{aligned}
\right.
\end{equation}
where $t\in [0,T]$ with $T>0$ being the deterministic terminal time;
$(\Omega,\mathcal{F}, \mathbb F, P)$ is a filtered complete probability space with
$\mathbb F = (\mathcal F_t)_{0\le t\le T}$ being the natural filtration of
the standard $d$-dimensional Brownian motion
$W = (W_t)_{0\le t \le T}$;
$X_0\in\F_0$ is the initial condition for the
forward stochastic differential equation (SDE);
$\xi\in \F_T$ is the terminal condition for the backward stochastic differential equation~(BSDE);
$b: \Omega\times [0,T]\times\R^q\times\R^p\times\R^{p\times d}\rightarrow\R^q$ and
$\sigma: \Omega\times [0,T]\times\R^q\times\R^p\times\R^{p\times d}\rightarrow\R^{q\times d}$ are
referred to the drift and diffusion coefficients respectively;
$f: \Omega\times [0.T]\times\R^q\times\R^p\times\R^{p\times d}\rightarrow\R^p$
is called the generator of BSDE and
$(X_t,Y_t,Z_t):[0,T]\times \Omega \rightarrow \R^q\times\R^p\times\R^{p\times d}$ is the unknown.
It is worth to note that $b(\cdot, x, y, z)$, $\sigma(\cdot,x, y, z)$ and $f(\cdot,x,y,z)$ are all
$\mathcal F_t$-adapted for any fixed numbers $x$, $y$ and $z$,
and that the two stochastic integrals with respect to $W_s$ are of the It\^o type. A triple $(X_t,Y_t,Z_t)$ is called an $L^2$-adapted solution
of the FBSDEs(\ref{FBSDEs}) if
it is $\F_t$-adapted, square integrable, and satisfies \eqref{FBSDEs}.
The FBSDEs \eqref{FBSDEs} is called \textit{decoupled} if $b$ and $\sigma$ are independent of $Y_t$ and $Z_t$.

Since the original work \cite{PP1}, in which Pardoux and Peng
obtained the existence and uniqueness of
the adapted solution for nonlinear BSDE under some standard conditions,
FBSDEs have been extensively studied. In \cite{MPJ},
a four-step approach was introduced to study the solvability
of coupled FBSDEs, and it was shown that the problem \eqref{FBSDEs} is uniquely solvable
under standard conditions on the data $b$, $\sigma$, $\xi$ and $f,$ when the diffusion
coefficient $\sigma$ does not depend on $Z_t$. Then Peng and Wu obtained the existence and
uniqueness of fully coupled FBSDE with an arbitrarily large time duration in \cite{PW}. 

In recent years, FBSDEs have possessed important applications
in many fields, such as mathematical finance \cite{KPQ,KKPPQ,MY1,KP}, risk measure \cite{P4, G1},
stochastic control \cite{P1,KKPPQ}, stochastic differential games \cite{HL}
and nonlinear expectations \cite{P3,G1}.
Generally speaking, it is difficult to find
the explicit closed-form solutions of FBSDEs, and thus, finding numerical solutions of FBSDEs becomes popular.
In the past decades, many numerical schemes for solving BSDE and decoupled FBSDEs have been proposed and studied.
Some of them, such as ~\cite{B1,BD,BP,BZ,BT,CZ,DM1,DM2,DMP,LGW,LGW1,MPMT,MSZ,MT,Z1},
are Euler-type methods with convergence rate $\frac{1}{2}$;
and some of them  \cite{cri,jfc1,jfc2,ZCP,ZLJ,ZLZ,ZWP,ZZJ,ZZwJ,SZZ} are high-order numerical methods with higher convergence rates. It is worth to point out that in the above literatures, the high-order methods for decoupled FBSDEs rely on the high-order approaches for both the forward SDE
and the backward SDE. It is clear that the high-order approaches for forward
SDEs requires large amounts of computations and are often difficult to be applied.
Concerning the coupled FBSDEs, since the coefficients of the forward SDE
depend on the unknowns $(Y_t, Z_t)$ of the BSDE, it seems not easy to design high-order (yet efficient) numerical schemes.

It is also worth noting that, in many practical problems such as the option pricing in stock markets,
one is only interested in the solutions of the BSDE, i.e. $Y_t$ and $Z_t$ in \eqref{FBSDEs}.
In fact, in the European option pricing problem,
the forward SDE is a mathematical model of stock prices, and its solution $X_t$ stands for
the price of the underlying assets at time $t,$ whose present value is regarded as a known quantity.
While $Y_t$ and $Z_t$ represent the option price and the
portfolio respectively. Therefore, what really make sense to us are the values of $Y_t$ and $Z_t$
at the known stock prices at the present time.
Such facts motivate us to use some \textit{cheap} numerical methods (such us the Euler method) to solve the forward SDE. Therefore, we are interested in the following problem:

\begin{itemize}
\item[\#] Whether one could still expect high-order accurate numerical solutions of the \textit{backward} SDE if the Euler method is used to discretize the \textit{forward} SDE?
\end{itemize}

In this paper, we shall give a positive answer to this question.
To this end, we shall study high accurate numerical methods for the coupled Markovian FBSDEs \eqref{FBSDEs},
in which the terminal condition $\xi$ is a deterministic
function $\varphi$ of $X_T$. Based on FBSDEs and SDE theories,
we will derive two reference ordinary differential equations (ODEs) of \eqref{FBSDEs} that contain the conditional expectations and their derivatives.
Then, our numerical scheme is obtained by approximating the conditional mathematical expectations
and their derivatives, in the reference ODEs. In particular, motivated by the local property of the generator of diffusion processes, the Euler method is used to solve the associated SDE in \eqref{FBSDEs}, which dramatically simplifies the entire computing complexity while without reducing the accuracy of the numerical solutions of the BSDE.
High-order convergence rates of such approaches are shown
by several numerical examples.

The rest of the paper is organized as follows. In Section \ref{pre},
some preliminaries on SDEs, FBSDEs and derivative approximations
are introduced. In Section \ref{scheme}, we discuss the design of high-order numerical methods for decoupled FBSDEs, and we extend these numerical methods to solve the coupled FBSDEs in Section \ref{scheme_c}. In Section \ref{NuEx}, numerical experiments are carried out to demonstrate the effectiveness of the proposed numerical schemes. We finally give some conclusions in Section \ref{cons}.

Some notations to be used:
\begin{itemize}
\item $\abs{\cdot}$: the Euclidean norm in the Euclidean space $\R$, $\R^q$ and $\R^{q\times d}$.
\item $\mathcal F_s^{t,x}$: the $\sigma$-algebra generated by the diffusion process
      $\{X_r,t\leq r\leq s, X_t=x\}$.
\item $\Ec{s}{t,x}{\eta}$: the conditional expectation of the random variable $\eta$
      under $\mathcal F_s^{t,x}$, i.e.,
      $\Ec{s}{t,x}{\eta} = \CE{\eta}{\mathcal F_s^{t,x}},$ and we use $\Ec{t}{x}{\eta}$ to denote $\Ec{t}{t,x}{\eta}$
      for simplicity.
\item $C_b^k:$ the set of functions $\phi(x)$ with uniformly bounded derivatives up
      to the order $k$.
\item $C^{k_1,k_2}:$ the set of functions $f(t,x)$ with continuous partial
derivatives up to $k_1$ with respect to $t$, and up to $k_2$ with respect to $x$.

\end{itemize}

\section{Preliminaries}\label{pre}


\subsection{Diffusion process and its generator}
A stochastic process $X_t$ is called a diffusion process starting at $(t_0, x_0)$ if it satisfies the SDE
\begin{equation}\label{SDEs}
\begin{aligned}
X_t =& x_0 + \int_{t_0}^t b(s,X_s) \di s + \int_{t_0}^t \sigma(s,X_s) \di W_s, \quad t\in[t_0, T],
\end{aligned}
\end{equation}
where $b_s=b(s,X_s)$ and $\sigma_s=\sigma(s,X_s)$ are measurable functions that satisfy
\begin{equation}\label{lcsde}\begin{aligned}
\abs{b(s,x)}+\abs{\sigma(s,x)} &\leq C(1+\abs{x}), \quad x\in \R^q, s\in [t_0,T],\\
\abs{b(t,x)-b(t,y)} + \abs{\sigma(t,x)-\sigma(t,y)} &\leq L\abs{x-y}, \quad x,y\in \R^q, s\in [t_0,T].
\end{aligned}\end{equation}
It is well known that under conditions \eqref{lcsde}, the SDE \eqref{SDEs}
has a unique solution. Moreover, one has $\Ec{t}{x}{X_s} = \CE{X_s}{X_t=x}, \forall t\leq s$,
by the Markov property of the diffusion process.

For a given measurable function $g: [0,T]\times \R^q\rightarrow \R$,
$g(t,X_t)$ is a stochastic process.
Let $G=\Ec{t}{x}{g(s,X_s)}$, then $G$ is a function of $(s,t,x)$ for $s\ge t$.
The partial derivative of $G$ with respect to $s$ is
defined by
\begin{equation}
\frac{\partial G}{\partial s} = \frac{\partial \Ec{t}{x}{g(s,X_s)}}{\partial s}
= \lim\limits_{\tau\downarrow 0} \frac{\Ec{t}{x}{g(s+\tau,X_{s+\tau})}-\Ec{t}{x}{g(s,X_s)}}{\tau},
\end{equation}
if the limit exists and is finite.

\begin{definition}\label{generator}
Let $X_t$ be a diffusion process in $\R^q$ satisfying \eqref{SDEs}.
The generator $A_t^x$ of $X_t$ on $g$ is defined by
\begin{equation}\label{gnrt}
A_t^x g(t,x) = \lim_{s\downarrow t}\frac{\Ec{t}{x}{g(s,X_s)}-g(t,x)}{s-t},\quad x\in\R^n.
\end{equation}
\end{definition}
%
%
Concerning the generator $A_t^x,$ we have the following result~\cite{O}:
\begin{thm}\label{gnrt_thm}
Let $X_t$ be the diffusion process defined by the SDE (\ref{SDEs}).
If $f\in C^{1,2}([0,T]\times\R^q)$, then we have
\begin{equation}\label{gnrt_dif}
A_t^x f(t,x) = L_{t,x}^0 f(t,x), \quad A_t^{X_t} f(t,X_t) = L_{t,X_t}^0 f(t,X_t),
\end{equation}
where \begin{equation}
L_{t,x}^0 = \frac{\partial}{\partial t} + \sum_i b_i(t,x)\frac{\partial}{\partial x_i}
+\frac12\sum_{i,j}\left(\sigma\sigma^\tau\right)_{i,j}(t,x)\frac{\partial^2}{\partial x_i\partial x_j}.
\end{equation}
\end{thm}
Note that $A_t^{X_t} f(t,X_t)\in \mathcal F_t$ is a stochastic process. Furthermore, by using together the It\^o's formula, Theorem \ref{gnrt_thm}
and the tower rule of conditional
expectations, we have the following theorem.


%

\begin{thm}\label{thm3}
Let $t_0<t$ be a fixed time, and $x_0\in \R^q$ be a fixed space point. If $f\in C^{1,2}([0,T]\times\R^q)$
and $\Ec{t_0}{x_0}{\abs{L_{t,X_t}^0f(t,X_t)}}<+\infty$,
we have
\begin{equation}\label{gnrt_eq}
\frac{\di\Ec{t_0}{x_0}{f(t,X_t)}}{\di t} = \Ec{t_0}{x_0}{A_t^{X_t}f(t,X_t)}, \quad t\ge t_0.
\end{equation}
Moreover, the following identity holds
\begin{equation}\label{gnrt_eq1}
\left. \frac{\di\Ec{t_0}{x_0}{f(t,X_t)}}{\di t} \right|_{t=t_0} = \left.\frac{\di\Ec{t_0}{x_0}{f(t,\bar X_t)}}{\di t}\right|_{t=t_0},
\end{equation}
where $\bar X_t$ is a diffusion process satisfying
\begin{equation}
\bar{X}_t = x + \int_{t_0}^t\bar{b}_s\di s + \int_{t_0}^t\bar\sigma_s\di W_s,
\end{equation}
and
$\bar b_s=\bar b(s,\bar X_s; t_0,x_0),$ $\bar \sigma_s = \bar{\sigma}(s,\bar X_s; t_0,x_0)$
are smooth functions of $(s,\bar X_s)$ with parameters $(t_0,x_0)$ that satisfy
\begin{equation*}
\bar b(t_0,\bar X_{t_0}; t_0,x_0)= b(t_0,x_0),\quad \bar{\sigma}(t_0,\bar X_{t_0}; t_0,x_0)=\sigma(t_0,x_0).
\end{equation*}
\end{thm}

Note that by choosing different $\bar b_s$ and $\bar\sigma_s,$
the identity \eqref{gnrt_eq1} gives different ways to approximate
$\left. \frac{\di\Ec{t_0}{x}{f(t,X_t)}}{\di t} \right|_{t=t_0}$.
The computational complexity can be reduced significantly if one uses
some appropriate choices. For example, one can simply choose
$\bar b(s,\bar X_s; t_0,x_0)= b(t_0,x_0)$ and
$\bar{\sigma}(s,\bar X_s; t_0,x_0)=\sigma(t_0,x_0)$ for all $s\in[t_0, t]$.




\subsection{Solution regularity and representation of FBSDEs}

Consider the following decoupled FBSDEs,

\begin{equation}\label{decoupled_fbsde}
\left\{\begin{aligned}
X_s^{t,x} &= x + \int_t^s b(r, X_r^{t,x}) \di r + \int_t^s \sigma(t, X_r^{t,x}) \di W_r,\\
Y_s^{t,x} &= \varphi(X_T^{t,x}) + \int_s^T f(r, X_r^{t,x}, Y_r^{t,x}, Z_r^{t,x})\di r
-\int_s^T Z_r^{t,x}\di W_r,
\end{aligned}\right. \forall s\in[t, T].
\end{equation}
Here the superscript $^{t,x}$ indicates that the forward SDE starts from $(t, x),$
which will be omitted if no ambiguity arises. To study the fine properties of decoupled
FBSDEs, the following assumptions were made in \cite{Z2}:
\begin{Assum}\label{assum1}
\begin{itemize}
\item[1.] The functions $b, \sigma \in C_b^1$, and there exists a constant $K > 0$,
such that
\[
\sup_{0\leq t\leq T}\{|b(t,0)| + |\sigma(t,0)|\} \leq K
\].
\item[2.] $q = d$, and the diffusion coefficient $\sigma$ is uniformly elliptic, that is,
there is a positive constant $K$ such that
\begin{equation}\label{sigma_ass}
\sigma(t, x)\sigma^T(t, x) \ge \frac{1}{K} I_q,\;\; \forall (t, x) \in [0; T] \times \R^d.
\end{equation}
\item[3.] The functions $f$ and $\varphi$ are uniformly Lipschitz continuous with Lipschitz
constant $K$, and assume
\[
\sup_{0\leq t\leq T}|f(t, 0, 0, 0)| + |\varphi(0)|\leq K
\]
\end{itemize}
\end{Assum}
Then, the following theorem was proved \cite{Z2}:
\begin{thm}\label{thm1} Under Assumption \ref{assum1}, there exists functions $u$ and $v$ and such that,
$\forall (t, x) \in [0, T) \times \R^d$, $Y_s^{t,x} = u(s,X_s^{t,x})$ and
\begin{enumerate}
\item $|v(t, x)|\le C$, where $C$ depends on $T$ and $K$;

\item $v$ is continuous;

\item $Z_s^{t,x} = v(s,X_s^{t,x}) \sigma(s,X_s^{t,x})$;

\item $\nabla_x u(t, x) = v(t, x)$;

\item If we assume further that $\varphi \in C_b^1$ , then 1--4 hold true on $[0, T] \times \R^d$, and
$v(T, x) = \nabla_x \varphi(x)$.
\end{enumerate}
\end{thm}

Moreover, we have the following nonlinear Feynman-Kac formula \cite{P2,MZ}:

\begin{thm}\label{thm2}  If the PDE
\begin{equation}\label{PDEs}
L_{t,x}^0u(t,x) + f(t,x,u(t,x),\nabla u(t,x)\sigma(t,x)) = 0
\end{equation}
with termial condition $u(T,x) = \varphi(x)$ has a classical solution
$u(t,x)\in C^{1,2}$, then
the unique solution $(X_s^{t,x},Y_s^{t,x},Z_s^{t,x})$ of the Markovian decoupled FBSDEs
\eqref{decoupled_fbsde}
with $\xi=\varphi(X_T^{t,x})$ can be represented as
\begin{equation}\label{F_k}
Y_s^{t,x} = u(s,X_s^{t,x}),\quad Z_s^{t,x} = \nabla_x u(s,X_s^{t,x})\sigma(s,X_s^{t,x}),\quad\forall s\in[t,T),
\end{equation}
where $\nabla_x u$ denotes the gradient of $u$ with respect to the spacial variable $x$.
\end{thm}


\subsection{Derivative approximation}\label{app_diff}
In this subsection, we will introduce a numerical method for approximating
the function derivatives (e.g., $\frac{d u(t)}{dt}$).
Such a method will play an
important role in designing our high-order numerical methods for FBSDEs.

Let $u(t)\in C_b^{k+1}$ and $t_i\in \mathbb R$ $(0\le i\le k)$ satisfying $t_0<t_1<\cdots<t_k$,
where $k$ is a positive integer. Let $\Delta  t_i=t_i-t_0,i=0,1,\ldots,k$.
Then by Taylor's expansion, for each $t_i,i = 0,1,\ldots,k$,
we have
\begin{flalign*}
\begin{split}
u(t_i) &= \sum_{j=0}^k\frac{(\Delta  t_i)^j}{j!}\frac{d^j u}{dt^j}(t_0) + \mathcal{O}\left(\Delta  t_i\right)^{k+1},
\end{split}
\end{flalign*}
then we can deduce
\begin{flalign*}
\begin{split}
\sum_{i=0}^k \alpha_{k,i} u(t_i) &= \sum_{j=0}^k\frac{\sum\limits_{i=0}^k\alpha_{k,i}(\Delta t_i)^j}{j!}\frac{d^j u}{dt^j}(t_0)
+\mathcal{O}\left(\sum_{i=0}^k\alpha_{k,i} (\Delta t_i)^{k+1}\right),
\end{split}
\end{flalign*}
where $\alpha_{k,i}$, $i=0,1,\dots,k$, are real numbers.
By choosing $\alpha_{k,i}~(i=0,1,\dots,k)$ as
\begin{equation}\label{alpha}
\frac{\sum\limits_{i=0}^k\alpha_{k,i} (\Delta  t_i)^j}{j!}
= \delta_{j1}=
\left\{
\begin{aligned}
& 1,\quad j=1,\\
& 0,\quad j\neq 1
\end{aligned}
\right.
\end{equation}
we obtain
\begin{equation}\label{dire_appro}
\frac{d u}{dt}(t_0) = \sum_{i=0}^k\alpha_{k,i} u(t_i)
+ Err,
\end{equation}
where
$Err=\mathcal{O}(\sum_{i=0}^k\alpha_{k,i} (\Delta t_i)^{k+1})$.
In particular, when $\Delta  t_i = i\Delta t$, the conditions (\ref{alpha})
are equivalent to the following linear system.
\begin{equation}\label{alpha_equ}
\left[
\begin{array}{ccccc}
1&1&1&\cdots&1\\
0&1&2&\cdots&k\\
0&1^2&2^2&\cdots&k^2\\
\vdots&\vdots&\vdots&\vdots&\vdots\\
0&1^k&2^k&\cdots&k^k
\end{array}
\right]\times
\left[
\begin{array}{c}
\alpha_{k,0}\Delta t\\
\alpha_{k,1}\Delta t\\
\alpha_{k,2}\Delta t\\
\vdots\\
\alpha_{k,k}\Delta t
\end{array}
\right]=
\left[
\begin{array}{c}
0\\
1\\
0\\
\vdots\\
0
\end{array}
\right].
\end{equation}
In the following table, 
we list $\alpha_{k,i}\Delta t$ $(i=0,1,\dots,k)$
of the system for $k=1,2,\ldots,6$.

\begin{table}[H]\label{alpha_value}
\setlength{\belowcaptionskip}{0pt}
\caption{}
\footnotesize
\begin{center}
\begin{tabular}{c|c|c|c|c|c|c|c}
\hline
$\alpha_{k,i}\Delta t $ & $i=0$ & $i=1$ & $i=2$ & $i=3$ & $i=4$ & $i=5$ & $i=6$\\
\hline
$k=1$ & $-1$              & $1$ &                 &                &                 &           & \\
\hline
$k=2$ & $-\frac32$        & $2$ & $-\frac12$      &                &                 &           & \\
\hline
$k=3$ & $-\frac{11}{6}$   & $3$ & $-\frac32$      & $\frac13$      &                 &           & \\
\hline
$k=4$ & $-\frac{25}{12}$  & $4$ & $-3$            & $\frac43$      & $-\frac14$      &           & \\
\hline
$k=5$ & $-\frac{137}{60}$ & $5$ & $-5$            & $\frac{10}{3}$ & $-\frac54$      & $\frac15$ & \\
\hline
$k=6$ & $-\frac{49}{20}$  & $6$ & $-\frac{15}{2}$ & $\frac{20}{3}$ & $-\frac{15}{4}$ & $\frac65$ & $-\frac16$\\
\hline
\end{tabular}\end{center}
\end{table}

\begin{remark}\label{rem:derivapp}
The multistep schemes proposed in this paper are closely related to the above
derivative approximation schemes. We now provide some remarks for the stability of the above numerical schemes, which have been well investigated in \cite{B2} for solving ODEs. To illustrate the
basic stability results, let us consider the ODE
\begin{equation}\label{ode:e1}
\frac{d  Y(t)}{d t} = f(t,Y(t)), \quad t\in[0,T)
\end{equation}
with the known terminal condition $Y(T)$.
Under the uniform time partition $0=t_0<t_1<\ldots<t_N=T$,
we apply the following multistep scheme to numerically solve \eqref{ode:e1}.
\begin{equation}\label{ode:e2}
\alpha_{k,0}Y^n + \sum_{j=1}^k\alpha_{k,j} Y^{n+j}=f(t_n,Y^n),
\end{equation}
where the $\alpha_{k,j}$'s are defined by \eqref {alpha_equ}.
By the theory of multistep scheme for solving ODEs \cite{B2},
in order to guarantee the stability of scheme \eqref{ode:e2},
the roots $\{\lambda_{k,j}\}_{j=1}^k$ of the charactristic equation
\begin{equation}\label{ode:e3}
P(\lambda)=\alpha_{k,0}\lambda^k + \sum_{j=1}^k\lambda^{k-j} =0
\end{equation}
must satisfy the root conditions, that is,
\[
|\lambda_{k,j}|\le 1.0, \,\, \textmd{and if} \,\, |\lambda_{k,j}|= 1.0, \,\, then \,\, P'(\lambda_{k,j})\ne 0
\,(\lambda_{k,j}\text{ are simple roots}).
\]
By the definition of $\alpha_{k,j}$ in \eqref{alpha_equ}, it can be checked
that $1$ is the simple root of \eqref{ode:e3} for each $k$.
In Table \ref{max_root1}, we list the maximum absolute values
of the roots for $k=2,3,\dots,8$ except the common root $1$.
\begin{table}[H]\label{max_root1}
\small
\setlength{\belowcaptionskip}{0pt}
\caption{The maximum absolute root of \eqref{ode:e3} except $1.0$}
\begin{center}
\begin{tabular}{c|c|c|c|c|c|c|c}
\hline
$k$ & 2 & 3 & 4 & 5 & 6 & 7 & 8 \\
\hline
$\max(\abs{\lambda_{k,j}})$  & 0.3333 & 0.4264 & 0.5608 & 0.7087 & 0.8633 & 1.0222 & 1.1839 \\
\hline
\end{tabular}
\end{center}
\end{table}
We learn from Table \ref{max_root1} that the multistep scheme \eqref{ode:e2}
is unstable for $k\geq 7,$
that is why we have only listed the $\alpha_{k,i}\Delta t$'s for $1\le k\le 6$ in Table \ref{alpha_value}.

\end{remark}



\section{Numerical schemes for decoupled FBSDEs}\label{scheme}
We first consider the numerical approaches for decoupled FBSDEs \eqref{FBSDEs},
namely, the functions $b$ and $\sigma$ are independent of $Y_t$ and $Z_t$.
Let $N$ be a positive integer.
For the time interval $[0,T]$, we introduce a regular time partition as
\[
0=t_0<t_1<t_2<\cdots<t_N=T.
\]
We will denote $t_{n+k}-t_n$ by $\Delta t_{n,k}$ and
$W_{t_{n+k}}-W_{t_n}$ by $\Delta W_{n,k},$
and use the notaions $\Delta t_{t_n,t}=t-t_n$ and $\Delta W_{t_n,t}=W_t-W_{t_n}$ for $t\geq t_n$.


%

\subsection{Two reference ODEs}

Let $(X_t,Y_t,Z_t)$ be the solution of the decoupled FBSDEs (\ref{FBSDEs})
with terminal condition $\xi = \varphi(X_T)$.
By taking conditional expectation $\Ec{t_n}{x}{\cdot}$ on both sides of the BSDE in (\ref{FBSDEs}),
we obtain the integral equation
\begin{equation}\label{exp_y}
\Ec{t_n}{x}{Y_t} = \Ec{t_n}{x}{\xi} + \int_t^T \Ec{t_n}{x}{f(s,X_s,Y_s,Z_s)}\di s,
\quad \forall t\in[t_n,T].
\end{equation}
By Theorem \ref{thm1},
the integrand $\Ec{t_n}{x}{f(s,X_s,Y_s,Z_s)}$ is a continuous
function of $s$.
Then, by taking derivative with respect to $t$ on both sides of \eqref{exp_y},
we obtain the following reference ordinary differential equation
\begin{equation}\label{dif_y}
\frac{\di \Ec{t_n}{x}{Y_t}}{\di t} = -\Ec{t_n}{x}{f(t,X_t,Y_t,Z_t)},\quad t\in[t_n,T].
\end{equation}

Note that we also have
$$Y_{t_n} = Y_t + \int_{t_n}^t f(s,X_s,Y_s,Z_s) \di s - \int_{t_n}^t Z_s \di W_s,\quad t\in[t_n,T].$$
By multiplying both sides of the above equation
by $(\Delta W_{t_n,t})^\tau$ (where $(\cdot )^\tau$
is the transpose of $(\cdot )$),
and taking the conditional expectation $\Ec{t_n}{x}{\cdot}$ on both sides of the derived equation,
we obtain
\begin{flalign}\label{exp_z}
\begin{split}
0 =~& \Ec{t_n}{x}{Y_{t}(\Delta W_{t_n,t})^\tau}
+\int_{t_n}^{t}\Ec{t_n}{x}{f(s,X_s,Y_s,Z_s)(\Delta W_{t_n,s})^\tau}\di s\\
&- \int_{t_n}^{t}\Ec{t_n}{x}{Z_s}\di s, \quad t\in[t_n,T].
\end{split}
\end{flalign}
Again,  the two integrands
in \eqref{exp_z} are continuous functions of $s$ by Theorem \ref{thm1}.
Upon taking derivative with respect to $t\in[t_n,T)$ in (\ref{exp_z})
one gets the following reference ODE:
\begin{equation}\label{dif_z}
\frac{\di \Ec{t_n}{x}{Y_t(\Delta W_{t_n,t})^\tau}}{\di t} = -\Ec{t_n}{x}{f(t,X_t,Y_t,Z_t)(\Delta W_{t_n,t})^\tau}
+ \Ec{t_n}{x}{Z_t}.
\end{equation}

\begin{remark}
The two ODEs \eqref{dif_y} and \eqref{dif_z} are our reference equations for the BSDE in (\ref{FBSDEs}). Our numerical schemes will be derived by approximating the conditional expectations and the derivatives in \eqref{dif_y} and \eqref{dif_z}.
\end{remark}

\subsection{The time semi-discrete scheme}
Motivated by Theorem \ref{thm3}, we choose smooth functions $\bar b(t,x)$
and $\bar\sigma(t,x)$ for $t\in [t_n,T]$ and $x\in \mathbb R^q$
with constraints $\bar b(t_n,x)=b(t_n,x)$ and $\bar \sigma(t_n,x)=\sigma(t_n,x)$.
The diffusion process $\bar X_t^{t_n,x}$ is defined by
\begin{equation}\label{xbar}\begin{aligned}
\bar{X}_t^{t_n,x} & = x + \int_{t_n}^t \bar b(s,\bar X_s^{t_n,x})\di s + \int_{t_n}^t \bar\sigma(s,\bar X_s^{t_n,x})\di W_s.
\end{aligned}\end{equation}

Let $(X_t^{t_n,x},Y_t^{t_n,x},Z_t^{t_n,x})$ be the solution of the decoupled FBSDEs.
According to Theorem \ref{thm2} and Theorem \ref{thm1},
the solutions $Y_t^{t_n,x}$ and $Z_t^{t_n,x}$
can be represented as $Y_t^{t_n,x}=u(t, X_t^{t_n,x})$ and
$Z_t^{t_n,x}=\nabla_x u(t, X_t^{t_n,x})\sigma(t,X_t^{t_n,x})$, respectively.

Let $\bar{Y}_t^{t_n,x}=u(t, \bar{X}_t^{t_n,x})$ and $\bar{Z}_t^{t_n,x}= \nabla_x u(t, \bar X_t^{t_n,x})\sigma(t,\bar X_t^{t_n,x}).$
By Theorem \ref{thm3}, we have
\begin{equation}\begin{aligned}
\left. \frac{\di\Ec{t_n}{x}{Y_t^{t_n,x}}}{\di t} \right|_{t=t_n}
&= \left.\frac{\di\Ec{t_n}{x}{\bar Y_t^{t_n,x}}}{\di t}\right|_{t=t_n},\\
\left. \frac{\di \Ec{t_n}{x}{Y_t^{t_n,x}(\Delta W_{t_n,t})^\tau}}{\di t}\right|_{t=t_n}
&= \left.\frac{\di \Ec{t_n}{x}{\bar Y_t^{t_n,x}(\Delta W_{t_n,t})^\tau}}{\di t}\right|_{t=t_n}.
\end{aligned}\end{equation}
Now introducing the scheme \eqref{dire_appro} into $\left.\frac{\di\Ec{t_n}{x}{\bar Y_t^{t_n,x}}}{\di t}\right|_{t=t_n}$ and
$\left.\frac{\di \Ec{t_n}{x}{\bar Y_t^{t_n,x}(\Delta W_{t_n,t})^\tau}}{\di t}\right|_{t=t_n}$,
we get

\begin{flalign}
\left.\frac{\di \Ec{t_n}{x}{Y_t^{t_n,x}}}{\di t}\right\vert_{t=t_n}
&= \sum_{i=0}^k\alpha_{k,i}\Ec{t_n}{x}{\bar Y_{t_{n+i}}^{t_n,x}} + \bar{R}_{y,n}^k,\label{app_y}\\
\left.\frac{\di \Ec{t_n}{x}{Y_t^{t_n,x}(\Delta W_{t_n,t})^\tau}}{\di t}\right\vert_{t=t_n}
&= \sum_{i=1}^k \alpha_{k,i} \Ec{t_n}{x}{\bar Y_{t_{n+i}}^{t_n,x}(\Delta W_{n,i})^\tau}
+\bar{R}_{z,n}^k, \label{app_z}
\end{flalign}
where $\alpha_{k,i}$ are defined by \eqref{alpha},
and $\bar{R}_{y,n}^k$ and $\bar{R}_{z,n}^k$ are truncation errors, i.e.
$$\bar{R}_{y,n}^k =\left.\frac{\di \Ec{t_n}{x}{Y_t^{t_n,x}}}{\di t}\right\vert_{t=t_n}
-\sum_{i=0}^k\alpha_{k,i}\Ec{t_n}{x}{\bar Y_{t_{n+i}}^{t_n,x}}, $$
$$\bar{R}_{z,n}^k = \left.\frac{\di \Ec{t_n}{x}{Y_t^{t_n,x}(\Delta W_{t_n,t})^\tau}}{\di t}\right\vert_{t=t_n}
- \sum_{i=1}^k \alpha_{k,i} \Ec{t_n}{x}{\bar Y_{t_{n+i}}^{t_n,x}(\Delta W_{n,i})^\tau}. $$

By inserting \eqref{app_y} and \eqref{app_z} into \eqref{dif_y} and \eqref{dif_z},
respectively, we obtain
\begin{flalign}
&\sum_{i=0}^k\alpha_{k,i}\Ec{t_n}{x}{\bar Y_{t_{n+i}}^{t_n,x}} = -f(t_n,x,Y_{t_n},Z_{t_n}) + R_{y,n}^{k},\label{y_ref}\\
&\sum_{i=1}^k\alpha_{k,i}\Ec{t_n}{x}{\bar Y_{t_{n+i}}^{t_n,x}(\Delta W_{n,i})^\tau} = Z_{t_n} + R_{z,n}^k,\label{z_ref}
\end{flalign}
where $R_{y,n}^k = -\bar R_{y,n}^k$ and $R_{z,n}^k= -\bar R_{z,n}^k$.

Let $Y^n$ and $Z^n$ be the numerical approximations
for the solutions $Y_t$ and $Z_t$ of the BSDE in \eqref{FBSDEs} at time $t_n$, respectively.
And we denoted by $X^n$ the numerical solution of the associated forward SDE at $t_n$.
Then, by removing the truncation error terms $R_{y,n}^k$ and $R_{z,n}^k$
from \eqref{y_ref} and \eqref{z_ref}, we get our time semi-discrete numerical scheme
for solving decoupled Markovian FBSDEs \eqref{FBSDEs}:

\begin{scheme}\label{semi_sch}
Assume that $Y^{N-i}$ and $Z^{N-i}$, $i = 0,1,\ldots,k-1$, are known.
For $n=N-k,\ldots,0$, with $\bar{X}_t^{t_n,X^n}$ being the solution of \eqref{xbar},
solve $Y^n=Y^n(X^n)$ and $Z^n=Z^n(X^n)$ by
\begin{flalign}
Z^n &= \sum_{j=1}^k\alpha_{k,j}\Ec{t_n}{X^n}{\bar Y^{n+j}(\Delta W_{n,j})^\tau},\label{semi_sch_z}\\
\alpha_{k,0}Y^n &= -\sum_{j=1}^k\alpha_{k,j}\Ec{t_n}{X^n}{\bar Y^{n+j}}
-f(t_n,X^n,Y^n,Z^n),\label{semi_sch_y}
\end{flalign}
where $\bar{Y}^{n+j}$ are the values of $Y^{n+j}$ at the space point $\bar{X}_t^{t_n,X^n}$.
\end{scheme}

The simplest choice of $\bar b$ and $\bar \sigma$ in \eqref{xbar} may be
that $\bar b(t,X_t^{t_n,x}) = b(t_n,x)$ and $\bar \sigma(t,X_t^{t_n,x}) = \sigma(t_n,x)$
for $t\in[t_n, T]$. In this case, Scheme \ref{semi_sch} becomes

\begin{scheme}\label{semi_sch_Eul}
Assume that $Y^{N-i}$ and $Z^{N-i}$, $i = 0,1,\ldots,k-1$, are known.
For $n=N-k,\ldots,0$, solve
$X^{n,j}(j=1,2,\dots,k)$, $Y^n=Y^n(X^n)$ and $Z^n=Z^n(X^n)$ by
\begin{flalign}
X^{n,j} = X^n + b(t_n, X^n)\Delta t_{n,j} + \sigma(t_n, X^n)\Delta W_{n,j},\quad j = 1,\ldots,k,\label{semi_sch_Eul_x}
\end{flalign}
\begin{flalign}
Z^n &= \sum_{j=1}^k\alpha_{k,j}\Ec{t_n}{X^n}{\bar Y^{n+j}(\Delta W_{n,j})^\tau},\label{semi_sch_Eul_z}\\
\alpha_{k,0}Y^n &= -\sum_{j=1}^k\alpha_{k,j}\Ec{t_n}{X^n}{\bar Y^{n+j}}
-f(t_n,X^n,Y^n,Z^n),\label{semi_sch_Eul_y}
\end{flalign}
where $\bar Y^{n+j}$ are the values of $Y^{n+j}$ at the space point $X^{n,j}$.
\end{scheme}

Note that by Theorem \ref{thm3}
and \eqref{dire_appro}, if $(L_{t,x}^0)^{k+1} u(t,x)$ is bounded,
we have \cite{B2}
\begin{equation}\label{trun_error_de}
\bar R_{y,n}^k=\mathcal{O}\left(\Delta t\right)^k,\quad \bar R_{z,n}^k = \mathcal{O}\left(\Delta t\right)^k,
\end{equation}
where $\bar R_{y,n}^k$ and $\bar R_{z,n}^k$ are defined
in \eqref{app_y} and \eqref{app_z}, respectively.


\begin{remark}
Motivated by Theorem 2, we used the Euler scheme to solve the associated SDE in Scheme
\ref{semi_sch_Eul}. The main advantages of this idea are twofold: for one hand, the use of the
Euler scheme can dramatically reduce the total computational complexity, and for the other hand,
one may be only interested in the solution of the BSDE in many applications.
We will show in Section \ref{scheme_c}
 that the numerical solutions of the BSDE can still be of
high-order accuracy in such a framework.
\end{remark}

\subsection{The fully discrete scheme}\label{Full_sch}
Scheme \ref{semi_sch_Eul} is a time semi-discrete scheme for solving FBSDEs.
To propose a fully discrete scheme, we introduce
a general space partition $D_h^n$ of $\R^q$ on each level $t_n$ with parameter $h^n>0$.
The space partition $D_h^n$ is a set of discrete grid points in $\R^q$, i.e
$D_h^n=\{x_i| x_i\in \R^q\}$.
We define the density of the grids in $D_h^n$ by
\begin{equation}\label{d_n^h}
h^n = \max\limits_{x\in\R^q}\min\limits_{x_i\in D_h^n} |x-x_i|=\max\limits_{x\in\R^q}\text{dist}(x,D_h^n),
\end{equation}
where $\text{dist}(x,D_h^n)$ is the distance from $x$ to $D_h^n$. For each $x \in \R^q,$ we define a local
subset $D^n_{h,x}$ of $D^n_h$ satisfying:
\begin{enumerate}
\item $\text{dist}(x,D_{h,x}^n) < \text{dist}(x,D_h^n\backslash D_{h,x}^n)$;
\item the number of elements in $D_{h,x}^n$ is finite and uniformly bounded, that is,
there exists a positive integer $\mathcal{N}_e$, such that, $\# D_{h,x}^n\leq \mathcal{N}_e$.
\end{enumerate}
We call $D_{h,x}^n$ the neighbor grid set in $D_h^n$ at $x$.

Now by \eqref{semi_sch_Eul_z} and \eqref{semi_sch_Eul_y}, we will solve $Y^n$ and $Z^n$
at grid points $x\in D_h^n$.
That is, for each $x\in D_{h}^n,n=N-k,\ldots,0$, we solve $Y^n=Y^n(x)$ and $Z^n=Z^n(x)$ by
\begin{flalign}
Z^n &= \sum_{j=1}^k\alpha_{k,j}\Ec{t_n}{x}{\bar Y^{n+j}(\Delta W_{n,j})^\tau},\label{semi_sch_z2}\\
\alpha_{k,0}Y^n &= -\sum_{j=1}^k\alpha_{k,j}\Ec{t_n}{x}{\bar Y^{n+j}}
-f(t_n,x,Y^n,Z^n),\label{semi_sch_y2}
\end{flalign}
where $\bar Y^{n+j}$ is the value of $Y^{n+j}$ at the space point $X^{n,j}$ defined by
\begin{equation}
X^{n,j} = X^n + b(t_n,X^n)\Delta t_{n,j} + \sigma(t_n,X^n)\Delta W_{n,j},\quad j = 1,\ldots,k.\label{semi_sch_x2}
\end{equation}
Generally, $X^{n,j}$ defined by  \eqref{semi_sch_x2}
does not belong to $D_{h}^{n+j}$ on condition of
 $X^n = x\in D_{h}^n$.
Thus, to solve $Y^n$ and $Z^n$, interpolation methods are needed to approximate
the value of $Y^{n+j}$ at $X^{n,j}$
using the values of $Y^{n+j}$ on $D_{h}^{n+j}$.
Here, we will adopt
a local interpolation operator $I_{h,X}^n$ such that $I_{h,X}^n g$ is the interpolation value of
the function $g$ at space point $X\in \R^q$ by using the values of $g$
only on $D_{h,X}^n$.
Note that any interpolation methods can be used here, however, care should be made if one wants to
guarantee the stability and accuracy.

In numerical simulations, the conditional expectations
$\Ec{t_n}{x}{\bar Y^{n+j}(\Delta W_{n,j})^\tau}$
and $\Ec{t_n}{x}{\bar Y^{n+j}}$ in \eqref{semi_sch_z2} and \eqref{semi_sch_y2}
should also be approximated.
The approximation operator of $\Ec{t_n}{x}{\cdot}$
is denoted by $\Ec{t_n}{x,h}{\cdot}$,
which can be any quadrature methods
such as the Monte-Carlo methods,
the quasi-Monte-Carlo methods, and the Gaussian quadrature methods and so on.

Now using the operators $I_{h,x}^n$ and $\Ec{t_n}{x,h}{\cdot}$,
we rewrite \eqref{y_ref} and \eqref{z_ref} in the following equivalent form.
\begin{equation}\label{FAP_YZ}\begin{array}{rl}
Z_{t_n} = & \sum_{j=1}^k\alpha_{k,j}\Ec{t_n}{x,h}{I_{h,\bar X_{t_{n+j}}}^{n+j} Y_{t_{n+j}}(\Delta W_{n,j})^\tau}\\
&-R_{z,n}^k + R_{z,n}^{k,\mathbb{E}} + R_{z,n}^{k,I_h},\\
\alpha_{k,0}Y_{t_n} =& -\sum_{j=1}^k\alpha_{k,j}\Ec{t_n}{x,h}{I_{h,\bar X_{t_{n+j}}}^{n+j} Y_{t_{n+j}}}
-f(t_n,x,Y_{t_n},Z_{t_n})\\
&+ R_{y,n}^k + R_{y,n}^{k,\mathbb{E}}+ R_{y,n}^{k,I_h},
\end{array}\end{equation}
where
$$\begin{array}{rl}
R_{z,n}^{k,\mathbb{E}} &= \sum_{j=1}^k\alpha_{k,j}(\mathbb E_{t_n}^{x}-\mathbb E_{t_n}^{x,h})\left[{\bar Y_{t_{n+j}}(\Delta W_{n,j})^\tau}\right],\\
R_{z,n}^{k,I_h} &= \sum_{j=1}^k\alpha_{k,j}\Ec{t_n}{x,h}{(\bar Y_{t_{n+j}}-I_{h,\bar X_{t_{n+j}}}^{n+j} Y_{t_{n+j}})(\Delta W_{n,j})^\tau},\\
R_{y,n}^{k,\mathbb{E}} &= -\sum_{j=1}^k\alpha_{k,j}(\mathbb E_{t_n}^{x}-\mathbb E_{t_n}^{x,h})\left[{\bar Y^{n+j}}\right],\\
R_{y,n}^{k,I_h} &= -\sum_{j=1}^k\alpha_{k,j}\Ec{t_n}{x,h}{\bar Y_{t_{n+j}}-I_{h,\bar X_{t_{n+j}}}^{n+j} Y_{t_{n+j}}}.
\end{array}$$
The two terms $R_{y,n}^{k,\mathbb{E}}$ and $R_{z,n}^{k,\mathbb{E}}$ are numerical errors
introduced by approximating conditional expectations,
and the other two terms $R_{y,n}^{k,I_h}$ and $R_{z,n}^{k,I_h}$ are numerical errors
caused by numerical interpolations.

By removing the six error terms $R_{y,n}^k$, $R_{z,n}^k$,
$R_{z,n}^{k,\mathbb{E}}$, $R_{z,n}^{k,I_h}$, $R_{y,n}^{k,\mathbb{E}}$
and $R_{y,n}^{k,I_h}$ from \eqref{FAP_YZ},
we obtain our fully discrete scheme for solving decoupled FBSDEs as follows:

\begin{scheme}\label{full_sch}
Assume $Y^{N-i}$ and $Z^{N-i}$ defined on $D_h^{N-i}$, $i = 0,1,\ldots,k-1$,
are known. For $n=N-k,\ldots,0$,
and for $x\in D_{h}^n$, solve $Y^n=Y^n(x)$ and $Z^n=Z^n(x)$ by
\begin{flalign}
X^{n,j} &= X^n + b(t_n,X^n)\Delta t_{n,j} + \sigma(t_n, X^n)\Delta W_{n,j},\quad j = 1,\ldots,k,\label{full_sch_x}\\
Z^n &= \sum_{j=1}^k\alpha_{k,j}\Ec{t_n}{x,h}{I_{h,X^{n,j}}^{n+j} Y^{n+j}(\Delta W_{n,j})^\tau},\label{full_sch_z}\\
\alpha_{k,0}Y^n &= -\sum_{j=1}^k\alpha_{k,j}\Ec{t_n}{x,h}{I_{h,X^{n,j}}^{n+j} Y^{n+j}}
-f(t_n,x,Y^n,Z^n).\label{full_sch_y}
\end{flalign}
\end{scheme}

Note that Scheme \ref{full_sch} is a $k$-step scheme.
For a fixed integer $k$, Scheme \ref{full_sch} consists of three procedures
for solving $Y^n$ and $Z^n$ at every $x\in D_h^n$ on each time level $t_n$:
(I) solve $X^{n,j}$ by the Euler scheme \eqref{full_sch_x};
(II) solve $Z^n$ by \eqref{full_sch_z} explicitly;
(III) solve $Y^n$ by \eqref{full_sch_y} implicitly.
Thus, some iterations are required for solving $Y^n.$
If the function $f(t_n,x,y,z)$
is Lipschitz continuous with respect to $y$, for small time partition step size $\Delta t_n$,
we can use the following
iteration procedure to approximately solve $Y^n$
\begin{flalign}
\alpha_{k,0}Y^{n,l+1} &= -\sum_{j=1}^k\alpha_{k,j}\Ec{t_n}{x,h}{I_{h,X^{n+j}}^{n+j} Y^{n+j}}
-f(t_n,x,Y^{n,l},Z^n),
\end{flalign}
until the iteration error $|Y^{n,l+1}-Y^{n,l}|\le \epsilon_0,$
where $\epsilon_0>0$ is a prescribed tolerance.
\begin{remark}\label{rem_trunc_err}
The local truncation errors of Scheme \ref{full_sch}
consist of six terms $R_{y,n}^k$, $R_{y,n}^{k,\mathbb{E}}$, $ R_{y,n}^{k,I_h}$,
$R_{z,n}^k$, $R_{z,n}^{k,\mathbb{E}}$ and $R_{z,n}^{k,I_h}$.
The two terms  $R_{y,n}^k$ and $R_{z,n}^k$ defined respectively in \eqref{y_ref} and \eqref{z_ref}
come from the approximations of the derivatives,
and the two terms $ R_{y,n}^{k,I_h}$ and $R_{z,n}^{k,I_h}$ defined in \eqref{FAP_YZ}
are the local interpolation errors.
For these terms,
suppose that the data $b$, $\sigma$, $f$ and $\varphi$
are sufficiently smooth such that $(L_{t,x}^0)^{k+1}u$ is bounded in $[0,T]\times \R^q$
and $u(t,\cdot)\in C_b^{r+1}, \forall t\in [0,T]$
and if the $k$-step scheme and $r$-degree polynomial interpolations are used,
then the following estimates hold \cite{AS,B2,BF}.
\begin{equation}\label{erro}
R_{y,n}^k=\mathcal{O}\left(\Delta t_n\right)^k, \,\, R_{z,n}^k=\mathcal{O}\left(\Delta t_n\right)^k, \,\,
R_{z,n}^{k,I_h}=\mathcal{O}\left(h^{r+1}\right), \,\, R_{y,n}^{k,I_h}=\mathcal{O}\left(h^{r+1}\right).
\end{equation}
The other two terms $R_{y,n}^{k,\mathbb{E}}$ and $R_{z,n}^{k,\mathbb{E}}$
are the local truncation errors resulted from the approximations of
the conditional mathematical expectations in \eqref{y_ref} and \eqref{z_ref}.
It is noted that
these conditional expectations
are functions of Gaussian random variables, which
can be represented as integrals with Gaussian kernels,
which may be approximated
by Gauss-Hermite Quadrature with high accuracy.
We will briefly introduce the Gauss-Hermite quadrature rule and its application
in Scheme \ref{full_sch} in the next subsection.
\end{remark}

In Scheme \ref{full_sch}, the values $\{(Y^{N-i},Z^{N-i})\}_{i=0}^{k-1}$ are needed
for solving $\{(Y^{n},Z^{n})\}^{N-k}_{n=0}.$ Note that $Y_T=\varphi(X_T),$ by Theorem \ref{thm1},
one has $Z_T=(\nabla_x \varphi \sigma)_{t=T}.$
Thus, we set $(Y^N,Z^N)=\big(\varphi(X_T), \nabla_x \varphi(X_T) \sigma(T,X_T)\big)$, however,
the values of $\{(Y^{N-i},Z^{N-i})\}_{i=1}^{k-1}$
are still needed to get an order-$k$ scheme.
In principle, these values can always obtained by solving the FBSDES in the interval
$[T-k\Delta t,T]$, with a \textit{more regular} time partition (fine mesh), as in \cite{ZZJ}.  Alternatively, one can sort to other types of high-order schemes, to get the values of $\{(Y^{N-i},Z^{N-i})\}_{i=1}^{k-1}.$ For instance, assume that $k=2,$ we can use
some simple schemes with \textit{local} truncation error of order $2,$ to get $(Y^{N-1},Z^{N-1}).$
For the cases of $k>2$, one can follow the similar procedure above.

\subsection{Gauss-Hermite quadrature rule for $\Ec{t_n}{x}{\cdot}$}\label{GHR}
The Gauss-Hermite quadrature rule is an extension of Gaussian quadrature
method for approximating the value of integrals of
$\int_{-\infty}^{+\infty}e^{-x^2}g(x)\di x$ by
\begin{equation}\label{ghf_1}
\int_{-\infty}^{+\infty}e^{-x^2}g(x)\di x \approx \sum_{j=1}^L\omega_jg(a_j),
\end{equation}
where $L$ is the number of sample points used in the approximation.
The points $\{a_j\}_{j=1}^L$ are the roots of the Hermite polynomial $H_L(x)$ of degree $L$
and $\{\omega_j\}_{j=1}^L$ are the corresponding weights \cite{AS}:
\[
\omega_j = \frac{2^{L+1}L!\sqrt{\pi}}{(H_L^\prime(a_j))^2}.
\]
The truncation error $R(g, L)$ of the Gauss-Hermite quadrature formula \eqref{ghf_1} is
\begin{equation}\label{trun_error_GHF}
R(g, L) = \int_{-\infty}^{+\infty}e^{-x^2}g(x)\di x - \sum_{j=1}^L\omega_jg(a_j)
 = \frac{L!\sqrt{\pi}}{2^L(2L)!}g^{(2L)}(\bar{x}),
\end{equation}
where $\bar{x}$ is a real number in $\R$. The Gauss-Hermite quadrature formula \eqref{ghf_1} is exact
for polynomial functions $g$ of degree less than $2L-1$.

For a $d$-dimensional function $g(\xx),\xx\in\R^d$, the Gauss-Hermite quadrature formula
becomes
\begin{equation}\label{ghf_q}
\int_{-\infty}^{\infty}\cdots \int_{-\infty}^{\infty}
g(\xx) e^{-\xx^\tau\xx}d\xx
\approx \sum_{\jj=\11}^{L}w_{\jj} g(\aa_{\jj}),
\end{equation}
where $\xx=(x_1,\dots,x_d)^\tau$, $\xx^\tau\xx=\sum\limits_{j=1}^d x_j^2$, and
$$
\mathbf{j} = (j_1,j_2,\ldots,j_d),\quad\omega_\mathbf{j} = \prod_{i = 1}^d\omega_{j_i},\quad \aa_{\mathbf{j}}=(a_{j_1},\ldots,a_{j_d})\quad\sum_{\mathbf{j} =1}^L =\sum_{j_1=1,\ldots,j_d=1}^{L,\ldots,L}.
$$
It is well known that, for a standard $d$-dimensional normal random variable $N(0,1)$, it holds that
\begin{equation}
\mathbb E[g(N)]=\frac{1}{(2\pi)^{\frac{d}{2}}}\int_{-\infty}^{+\infty}g(\xx)e^{-\frac{\xx^\tau \xx}{2}}\di \xx
=\frac{1}{(\pi)^{\frac{d}{2}}}\int_{-\infty}^{+\infty}g(\sqrt{2}\xx)e^{-\xx^\tau\xx}\di \xx.
\end{equation}
Then by \eqref{ghf_q}, we get
\begin{equation}\label{E_GH}
\mathbb E[g(N)] =  \frac{1}{(\pi)^{\frac{d}{2}}} \sum_{\jj=\11}^{L}w_{\jj} g(\aa_{\jj})+R_{\mathbb{E},L}^{GH}(g),
\end{equation}
where $R_{\mathbb E,L}^{GH}(g) = \mathbb E[g(N)] -  \frac{1}{(\pi)^{\frac{d}{2}}} \sum_{\jj=\11}^{L}w_{\jj} g(\aa_{\jj})$
is the truncation error of the Gauss-Hermite quadrature rule for $g$.

Recall that, in Scheme \ref{full_sch}, the conditional expectation
$\Ec{t_n}{\xx}{\bar Y^{n+j}}$
is approximated by $\Ec{t_n}{\xx,h}{I_{h,X^{n+j}}^{n+j} Y^{n+j}}$, where
$\Ec{t_n}{\xx,h}{\cdot}$ is the approximation of $\Ec{t_n}{\xx}{\cdot}$,
 and $I_{h,X^{n+j}}^{n+j} Y^{n+j}$ is the interpolation approximation of $\bar Y^{n+j}$.
By the nonlinear Feynman-Kac formula \eqref{F_k} and
Scheme \ref{full_sch},
$\bar Y^{n+j}$ is a function of $X^{n,j}$
and has the following explicit representation.
\begin{equation}
\bar Y^{n+j} = Y^{n+j}(X^{n+j}) = Y^{n+j}(X^n + b(t_n,X^n)\Delta t_{n,j} + \sigma(t_n, X^n)\Delta W_{n,j}),
\end{equation}
where $\Delta W_{n,j}\sim\sqrt{\Delta t_{n,j}}N(0,I_d)$
is a $d$-dimensional Gaussian random variable.
Thus we define the approximation
$\Ec{t_n}{\xx,h}{\bar Y^{n+j}}$ by
\begin{equation}\label{E_GHhy}
\Ec{t_n}{\xx,h}{\bar Y^{n+j}} =  \frac{1}{(\pi)^{\frac{d}{2}}}
\sum_{\jj=\11}^{L}w_{\jj} FY(\aa_{\jj}),
\end{equation}
where $FY=FY(y)=Y^{n+j}(\xx + b(t_n,\xx)\Delta t_{n,j} + \sigma(t_n, \xx)\sqrt{2\Delta t_{n,j}} y)$.
We denote such approximation error by $R_{\mathbb E_{t_n}^{\xx},L}^{GH}(FY).$

Similarly, we have
\begin{equation}\label{E_GHhz}
\Ec{t_n}{\xx,h}{\bar Y^{n+j}(\Delta W_{t_n,j})^\tau} =  \frac{1}{(\pi)^{\frac{d}{2}}}
\sum_{\jj=\11}^{L}w_{\jj} FYW(\aa_{\jj}),
\end{equation}
where $FYW=FYW(y)=Y^{n+j}(\xx + b(t_n,\xx)\Delta t_{n,j} + \sigma(t_n,\xx)\sqrt{2\Delta t_{n,j}} y)y$.
The approximation error is denoted by $R_{\mathbb E_{t_n}^{\xx},L}^{GH}(FYW).$

For smooth data, by \eqref{trun_error_GHF}, we have the following estimates \cite{AS,STW}
\begin{equation}\label{er_mc}
R_{y,n}^{k,\mathbb{E}}=\mathcal{O}\left(\frac{L!}{2^L(2L)!}\right),\quad R_{z,n}^{k,\mathbb{E}}=\mathcal{O}\left(\frac{L!}{2^L(2L)!}\right),
\end{equation}
where $R_{y,n}^{k,\mathbb{E}}$ and $R_{z,n}^{k,\mathbb{E}}$ are defined in \eqref{FAP_YZ}.

Now, we would like to remark that, to obtain a high-order scheme for solving FBSDEs, both the interpolation error and the above integration error should be well controlled, to balance the time discretization error.

\section{Numerical schemes for coupled FBSDEs}\label{scheme_c}

In this section we extend Scheme \ref{full_sch} to the following scheme for
solving fully coupled FBSDEs \eqref{FBSDEs}.

\begin{scheme}\label{full_nsch1}
Assume $Y^{N-i}$ and $Z^{N-i}$ defined on $D_h^{N-i}$, $i = 0,1,\ldots,k-1$,
are known. For $n=N-k,\ldots,0$,
and for $x\in D_{h}^n$, solve $Y^n=Y^n(x)$ and $Z^n=Z^n(x)$ by
\begin{flalign}
X^{n,j} =& X^n + b(t_n,X^n,Y^n,Z^n)\Delta t_{n,j} + \sigma(t_n,X^n,Y^n,Z^n)\Delta W_{n,j},\nonumber  \\
             &  j = 1,2,\ldots,k,\nonumber \\
Z^n = &\sum_{j=1}^k\alpha_{k,j}\Ec{t_n}{x,h}{I_{h,X^{n,j}}^{n+j} Y^{n+j}(\Delta W_{n,j})^\tau},\nonumber\\
\alpha_{k,0}Y^n = & -\sum_{j=1}^k\alpha_{k,j}\Ec{t_n}{x,h}{I_{h,X^{n,j}}^{n+j} Y^{n+j}}
-f(t_n,x,Y^n,Z^n).\nonumber
\end{flalign}
\end{scheme}

In Scheme \ref{full_nsch1}, $X^{n,j}$, $Y^n$ and $Z^n$ are
coupled together. Numerical methods for solving nonlinear
equations are needed. In our numerical experiments,
we will use the following iterative scheme to solve
$Y^n$ and $Z^n$.

\begin{scheme}\label{full_nsch}
Assume $Y^{N-i}$ and $Z^{N-i}$ defined on $D_h^{N-i}$, $i = 0,1,\ldots,k-1$,
are known. For $n=N-k,\ldots,0$,
and for $x\in D_{h}^n$, solve $Y^n=Y^n(x)$ and $Z^n=Z^n(x)$ by
\begin{enumerate}
\item let $Y^{n,0}=Y^{n+1}(x)$ and $Z^{n,0}=Z^{n+1}(x)$,
and let $l=0$;
\item  for $l=0,1,\dots$, solve $Y^{n,l+1}=Y^{n,l+1}(x)$ and $Z^{n,l+1}=Z^{n,l+1}(x)$ by
\begin{flalign}
X^{n,j} = & X^n + b(t_n,X^n,Y^{n,l},Z^{n,l})\Delta t_{n,j} + \sigma(t_n,X^n,Y^{n,l},Z^{n,l})\Delta W_{n,j}, \nonumber \\
& \,\, j = 1,2,\ldots,k,\nonumber \\
Z^{n,l+1} = & \sum_{j=1}^k\alpha_{k,j}\Ec{t_n}{x,h}{I_{h,X^{n,j}}^{n+j} Y^{n+j}(\Delta W_{n,j})^\tau},\nonumber\\
\alpha_{k,0}Y^{n,l+1} = & -\sum_{j=1}^k\alpha_{k,j}\Ec{t_n}{x,h}{I_{h,X^{n,j}}^{n+j} Y^{n+j}}
-f(t_n,x,Y^{n,l+1},Z^{n,l+1})\nonumber
\end{flalign}
until $\max(|Y^{n,l+1}-Y^{n,l}|,|Z^{n,l+1}-Z^{n,l}|)< \epsilon_0$;
\item Let $Y^n = Y^{n,l+1}$ and $Z^n=Z^{n,l+1}$.
\end{enumerate}
\end{scheme}

Note that if the drift coefficient $b$ and the diffusion
coefficient $\sigma$ do not depend on $Y_.$ and $Z_.$, Scheme \ref{full_nsch}
coincides with Scheme \ref{full_sch}.

\begin{remark}
Notice that the mesh $D^n_h$ is essentially unbounded.
However, in real computations, we are interested in obtaining only certain values of $(Y_t, Z_t)$
at $(t_n,x)$ with $x$ belongs to a bounded domain.
For instance,  in option pricing, we are only interested in the option values at the current
option price. Thus,  in practice,
only a bounded sub-mesh of $D^n_h$ is used on each time level.
It is also worth to note that, in our numerical tests, we use the Gaussian-Hermite integral rule to approximate
the conditional expectations. Such an approximation is global, while only a small number of integral points are needed (of course, the integral points are bounded).
\end{remark}

\section{Numerical experiments}\label{NuEx}
In this section, we will provide with several numerical examples to show the behavior of our
fully discrete Scheme \ref{full_sch} and Scheme \ref{full_nsch} for solving FBSDEs.
For simplicity, we will use uniform partitions in both time and space. The time interval $[0,T]$ will be uniformly divided into $N$ parts with time step $\Delta t=\frac{T}{N}$
and time grids $t_n = n\Delta t$, $n=0,1,\dots,N$. The space partition is $D_h^n=D_h$ for all $n$, where
$D_h$ is defined by
\begin{equation}
D_h = D_{1,h} \times D_{2,h} \times \cdots \times D_{q,h},
\end{equation}
where $D_{j,h}$ is the partition of the one-dimensional real axis $\R$
\begin{equation*}
D_{j,h} = \left\{ x_i^j: x_i^j = ih, i =0,\pm 1,\cdots, \pm\infty\right\}
\end{equation*}
for $j = 1,2,\ldots,q$, and $D_{h,x}\subset D_h$ denotes the set
of some neighbor grids near $x$. For all the examples, the terminal time $T$
is set to be $1.0.$

In our numerical experiments, the approximation of conditional expectations
$\Ec{t_n}{x,h}{\bar Y^{n+j}}$ and $\Ec{t_n}{x,h}{\bar Y^{n+j}\Delta W_{n,j}}$
were defined by \eqref{E_GHhy} and \eqref{E_GHhz} respectively for $x\in D_h$.
We choose $I_{h,x}^{n}$ to be the local Lagrange interpolation operator
based on the values of the interpolated functions
on $D_{h,x}$ at time level $t_n$,
so that the interpolation error estimates in \eqref{erro} hold.

Our main goal of the numerical experiments
is to demonstrate the high accuracy of
the fully discrete Scheme \ref{full_sch} and Scheme \ref{full_nsch}
in time.
So in the conditional expectation approximations,
we set the number of the Gauss-Hermite quadrature points to be
big enough (we take eight points in each dimension) such that the errors resulted from
the use of Gauss-Hermite quadrature rule are negligible.
To balance the time discrete truncation errors
and the space truncation error,
for fixed time partition step size $\Delta t$,
we require $h^{r+1}=(\Delta t)^{k+1}$, where $k$ is the step number of our scheme and $r$ is
the degree of the Lagrangian interpolation polynomials.
In our numerical experiments, we change the problem
of choosing $r$ into the problem of choosing $h$ such that $h = \Delta t^{\frac{k+1}{r+1}}$.
we choose small $r$ when
lower order schemes $(k\leq 3)$ are tested (such as $r=4$ or $r=6$). And when the order becomes higher
than $3$, we choose a bigger $r$ (e.g., $r=10$ or $r=15$).
We assume that $\{(Y^{N-j},Z^{N-j})\}_{j=1}^k$ are
given for fixed $k$
in such a way that their effects for the the convergence rate are also negligible. However, as discussed before, the value of $\{(Y^{N-j},Z^{N-j})\}_{j=1}^k$ can also be computed numerically.

In what follows, we denote by CR and RT the convergence rate and the running time
respectively. The unit of RT is second. The numerical results, including numerical errors, convergence rates and running times,
are obtained by running Scheme \ref{full_sch} and Scheme \ref{full_nsch} in FORTRAN 95
on a computer with 16 Intel Xeon E5620 CPUs(2.40GHz), and 3.1G free RAM.

\subsection{Decoupled Cases} \;\;

In this subsection we use two numerical examples
to demonstrate the high accuracy
of the fully discrete Scheme \ref{full_sch} for solving decoupled
FBSDEs.

\begin{example}\label{EX_NONL} {\rm
The considered decoupled FBSDEs are

\begin{equation}\label{NONL_fbsde}
\left\{
 \begin{array}{rl}
    \di X_t =& \frac{1}{1+2\exp(t+X_t)}\di t + \frac{\exp(t+X_t)}{1+\exp(t+X_t)}\di W_t,\\
    -\di Y_t =& \left(-\frac{2Y_t}{1+2\exp(t+X_t)}-\frac12\left(\frac{Y_tZ_t}{1+\exp(t+X_t)}-Y_t^2Z_t\right)\right)\di t
               -Z_t\di W_t
 \end{array}
\right.
\end{equation}
with the initial and terminal conditions
$X_0 = x$, $Y_T = \frac{\exp(t+X_T)}{1+\exp(t+X_T)}$.

The analytic solutions $Y_t$ and $Z_t$ of \eqref{NONL_fbsde} are
\begin{equation}
\left\{
 \begin{aligned}
    Y_t &= \frac{\exp(t+X_t)}{1+\exp(t+X_t)},\\
    Z_t &= \frac{(\exp(t+X_t))^2}{(1+\exp(t+X_t))^3}.
 \end{aligned}
\right.
\end{equation}

We choose the initial condition $x=1.0$,
and solve the approximate solutions $Y^0$ and $Z^0$ by Scheme \ref{full_sch}
for $k=1,2,\dots,8.$ The errors $|Y^0-Y_0|$
and $|Z^0-Z_0|$, the convergence rates (CR), and running times (RT)
are listed in Table \ref{tableEX_NONL-1}

\begin{table}[H]\label{tableEX_NONL-1}
\small
\setlength{\belowcaptionskip}{0pt}
\caption{Errors and convergence rates for Example \ref{EX_NONL}.}
\begin{center}
\begin{tabular}{c|c|c|c|c|c|c|c}
\hline
Scheme  \ref{full_sch}                &             &   N=16   &   N=32   &   N=64   &   N=128  &   N=256  & CR         \\
\hline
\multirow{3}{*}{$k=1$} & $|Y^0-Y_0|$ & 3.576E-03 & 1.789E-03 & 8.946E-04 & 4.474E-04 & 2.238E-04 & 1.000 \\
                        \cline{2-8}
                        & $|Z^0-Z_0|$ & 4.322E-03 & 2.161E-03 & 1.081E-03 & 5.405E-04 & 2.703E-04 & 1.000 \\
                        \cline{2-8}
                        & RT          & 0.236     & 0.537     & 1.155     & 4.090     & 23.479    &       \\
\hline
\multirow{3}{*}{$k=2$} & $|Y^0-Y_0|$ & 8.199E-05 & 2.140E-05 & 5.458E-06 & 1.378E-06 & 3.460E-07 & 1.973 \\
                        \cline{2-8}
                        & $|Z^0-Z_0|$ & 1.103E-04 & 2.678E-05 & 6.585E-06 & 1.631E-06 & 4.053E-07 & 2.021 \\
                        \cline{2-8}
                        & RT          & 0.443     & 0.628     & 1.871     & 9.631     & 62.101    &       \\
\hline
\multirow{3}{*}{$k=3$} & $|Y^0-Y_0|$ & 6.377E-07 & 8.299E-08 & 1.024E-08 & 1.264E-09 & 1.567E-10 & 3.002 \\
                        \cline{2-8}
                        & $|Z^0-Z_0|$ & 5.456E-06 & 8.108E-07 & 1.099E-07 & 1.424E-08 & 1.818E-09 & 2.893 \\
                        \cline{2-8}
                        & RT          & 0.617     & 1.011     & 3.310     & 18.367    & 144.191   &       \\
\hline
\multirow{3}{*}{$k=4$} & $|Y^0-Y_0|$ & 1.446E-07 & 1.026E-08 & 6.795E-10 & 4.388E-11 & 2.761E-12 & 3.922 \\
                        \cline{2-8}
                        & $|Z^0-Z_0|$ & 1.314E-06 & 9.316E-08 & 6.181E-09 & 3.996E-10 & 2.530E-11 & 3.919 \\
                        \cline{2-8}
                        & RT          & 0.768     & 1.467     & 5.705     & 35.749    & 264.044   &       \\
\hline
\multirow{3}{*}{$k=5$} & $|Y^0-Y_0|$ & 7.821E-08 & 2.379E-09 & 8.836E-11 & 2.725E-12 & 3.492E-14 & 5.196 \\
                        \cline{2-8}
                        & $|Z^0-Z_0|$ & 7.544E-07 & 2.482E-08 & 7.590E-10 & 2.358E-11 & 6.860E-13 & 5.017 \\
                        \cline{2-8}
                        & RT          &   1.744   &   5.005   &   21.005   &   95.826  &  679.963  &       \\

\hline
\multirow{3}{*}{$k=6$} & $|Y^0-Y_0|$ & 1.095E-08 & 1.827E-10 & 4.167E-12 & 1.960E-14 & 2.109E-14 & 5.116(6.273) \\
                        \cline{2-8}
                        & $|Z^0-Z_0|$ & 5.427E-08 & 5.464E-10 & 5.929E-11 & 6.004E-14 & 1.427E-14 & 5.687(6.256) \\
                        \cline{2-8}
                        & RT          & 2.049     & 3.492     & 17.068     & 113.369    &  907.924  &       \\
\hline
\multirow{3}{*}{$k=7$} & $|Y^0-Y_0|$ & 1.963E-08 & 2.532E-10 & 3.943E-12 & 6.506E-14 & 3.105E-13 & 4.382 \\
                        \cline{2-8}
                        & $|Z^0-Z_0|$ & 9.879E-08 & 2.713E-09 & 2.826E-11 & 4.254E-13 & 5.416E-13 & 4.759 \\
                        \cline{2-8}
                        & RT          &  3.183    &  8.299    &   45.778   &  289.968   &  2700.751  &       \\
\hline
\multirow{3}{*}{$k=8$} & $|Y^0-Y_0|$ & 6.886E-07 & 7.066E-09 & 3.141E-09 & 7.148E-04 & 3.923E-01 & -5.487 \\
                        \cline{2-8}
                        & $|Z^0-Z_0|$ & 8.200E-06 & 9.859E-08 & 5.371E-08 & 2.145E-02 & 1.092E+03 & -7.170 \\
                        \cline{2-8}
                        & RT          & 6.630     & 21.202     & 107.122     & 536.644    &  3284.037  &       \\
\hline
\end{tabular}
\end{center}
\end{table}
}
\end{example}
From Table \ref{tableEX_NONL-1} we can make the following conclusions:
\begin{enumerate}
\item Scheme \ref{full_sch} is a $k$-order numerical scheme
for $1\le k\le 6$ (note that when $k=6$, the convergence rate for
$N=16, 32, 64$ and $128$ is approximately $6.0$, but when $N=256$, the accuracy of
double precision variables influenced the convergence rate).
Such result is consistent with the one
when this kind of multistep scheme is used to solve deterministic ODEs \cite{B2}.

\item The errors and the convergence rates listed in Table \ref{tableEX_NONL-1} show
that Scheme \ref{full_sch} is stable for $1\le k\le 6$, and is unstable when $k\geq 7,$ which
is again consistent with those in deterministic ODEs settings \cite{B2}. We would like to mention that the $L^2$-stability for a large class of linear
multistep schemes has been studied in \cite{jfc1}, however, the analysis in \cite{jfc1} does not cover our new schemes. A rigorous stability proof for our schemes is part of our ongoing work.

\item The errors and the running times listed in Table \ref{tableEX_NONL-1} clearly show that
 Scheme \ref{full_sch} will be more efficient if bigger $k$ is used
 ($1\le k\le 6$) in general.
\end{enumerate}

Generally, a highly accurate numerical scheme
for solving the forward SDE is needed, if one wants to obtain a highly accurate numerical scheme
for the BSDE. However, the results in Table \ref{tableEX_NONL-1} indeed show that highly accurate numerical solutions of the BSDE are obtained when the
Euler scheme is used to solve the forward SDE. Such idea can
significantly simplify the solution procedure of the scheme for solving FBSDEs. To the best of our knowledge, no such effort has been made before for solving FBSDEs.

Let's now take the European option pricing as an example, for which one is only interested in the solution of the BSDE.
\begin{example}\label{EX_BS}{\rm Assume that the option price $S_t$ and the bond price $p_t$ are governed by
\begin{equation}\label{BS_stock}
\begin{aligned}
S_t &= S_0 + \int_0^t b_\tau S_\tau\di \tau + \int_0^t \sigma_\tau S_\tau\di W_\tau, \quad \tau\geq 0\\
p_t &= p + \int_0^t r_\tau p_\tau\di \tau, \quad \tau\geq 0,\\
\end{aligned}
\end{equation}
where $b_t$ is the expected return rate of the stock, $\sigma_t>0$ is the volatility of the
stock, $r_t$ is the return rate of the bond, and $b_t$, $r_t$, $\sigma_t$ and $\sigma_t^{-1}$
are all bounded.

Suppose that
an investor with wealth $Y_t$ puts $\pi_t$ money to buy the stock and use $Y_t-\pi_t$ to buy
the bond, and the stock pays dividends continuously with a bounded dividend rate $d(t,S_t)$
at time $t$. Then by using the no-arbitrage principle,
the processes $Y_t$ and $\pi_t$ satisfy the following stochastic
differential equation.
\begin{equation}
-\di Y_t = -(r_t Y_t + (b_t-r_t+d(t,S_t))\pi_t)\di t - \sigma_t\pi_t\di W_t.
\end{equation}
Let $Z_t = \sigma_t\pi_t$, then $(Y_t, Z_t)$ satisfies
\begin{equation}\label{BS_bsde}
-\di Y_t = -(r_t Y_t + (b_t-r_t+d(t,S_t))\sigma_t^{-1}Z_t)\di t - Z_t\di W_t.
\end{equation}
For the European call option the terminal condition is given at the mature time
$T$ by
\begin{equation}\label{BS_ter}
Y_T = \max\{S_T-K, 0\},
\end{equation}
where $S_T$ is the solution $S_t$ of the model \eqref{BS_stock} at the mature time $T$, and
$K$ is the strike price. Therefore, by the no-arbitrage principle, $Y_0$ should be the price
of this European call option.

In particular, when $b_t=b$, $\sigma_t=\sigma$, $r_t=r$ and $d(t, S_t)=d$ are all constants
in \eqref{BS_stock} and \eqref{BS_bsde}, by the Black-Scholes formula, the analytic solution
$(Y_t,Z_t)$ of the FBSDE \eqref{BS_stock} and \eqref{BS_bsde} with terminal condition \eqref{BS_ter}
can be written in the following form.
\begin{equation}\label{BS_true}
\left\{
\begin{array}{rl}
Y_t &= V(t, S_t) = S_t e^{-d(T-t)}N(d_1(S_t)) - K e^{-r(T-t)}N(d_0(S_t)),\\
Z_t &= \frac{\partial V}{\partial S}\sigma = S_t e^{-d(T-t)}N(d_1(S_t))\sigma,\\
d_0(S_t) &= \frac{1}{\sigma\sqrt{T-t}}\log\{\frac{S_t}{K e^{(d-r)(T-t)}}\}-\frac12\sigma\sqrt{T-t},\\
d_1(S_t) &= d_0(S_t) + \sigma\sqrt{T-t},
\end{array}
\right.
\end{equation}
where $N$ is the cumulative normal distribution function.

Now we take the constants $b=0.05$, $\sigma=0.2$, $r=0.03,$ $d=0.04$, the mature time $T=1.0$ and
$K=S_0=100.0$. We solve the problem by Scheme \ref{full_sch}
for $k=1,2,3,4$. The errors, convergence
rates and running time are listed in Table \ref{tableBS-1}.

\begin{table}[H]\label{tableBS-1}
\small
\setlength{\belowcaptionskip}{0pt}
\caption{Errors and convergence rates for Example \ref{EX_BS}.}
\begin{center}
\begin{tabular}{c|c|c|c|c|c|c|c}
\hline
Scheme                  &             &   N=16   &   N=32   &   N=64   &   N=128  &   N=256  & CR         \\
\hline
\multirow{3}{*}{$k=1$} & $|Y^0-Y_0|$ & 1.401E-02 & 6.988E-03 & 3.490E-03 & 1.744E-03 & 8.716E-04 & 1.002 \\
                        \cline{2-8}
                        & $|Z^0-Z_0|$ & 1.795E-01 & 8.951E-02 & 4.469E-02 & 2.233E-02 & 1.116E-02 & 1.002 \\
                        \cline{2-8}
                        & RT          & 0.204     & 0.412     & 1.160     & 5.554     & 37.673    &       \\
\hline
\multirow{3}{*}{$k=2$} & $|Y^0-Y_0|$ & 6.181E-04 & 1.637E-04 & 4.207E-05 & 1.066E-05 & 2.682E-06 & 1.964 \\
                        \cline{2-8}
                        & $|Z^0-Z_0|$ & 1.311E-03 & 3.286E-04 & 8.233E-05 & 2.061E-05 & 5.155E-06 & 1.998 \\
                        \cline{2-8}
                        & RT          & 0.374     & 0.744     & 2.504     & 13.224    & 93.017    &       \\
\hline
\multirow{3}{*}{$k=3$} & $|Y^0-Y_0|$ & 3.242E-05 & 4.502E-06 & 5.905E-07 & 7.554E-08 & 9.551E-09 & 2.935 \\
                        \cline{2-8}
                        & $|Z^0-Z_0|$ & 1.416E-05 & 2.429E-06 & 3.453E-07 & 4.582E-08 & 5.894E-09 & 2.819 \\
                        \cline{2-8}
                        & RT          & 0.536     & 1.244     & 4.160     & 25.070    & 217.832   &       \\
\hline
\multirow{3}{*}{$k=4$} & $|Y^0-Y_0|$ & 1.808E-06 & 1.755E-07 & 1.059E-08 & 6.061E-10 & 3.403E-11 & 3.957 \\
                        \cline{2-8}
                        & $|Z^0-Z_0|$ & 5.174E-06 & 4.030E-07 & 2.419E-08 & 1.460E-09 & 8.469E-11 & 3.991 \\
                        \cline{2-8}
                        & RT          & 0.764     & 1.895     & 6.966     & 48.928    & 367.599   &       \\
\hline
\end{tabular}
\end{center}
\end{table}
}
\end{example}

Again, it is shown in Table \ref{tableBS-1} that
Scheme \ref{full_sch} is a $k$-order scheme. Furthermore,
the scheme with a larger $k$ is more competitive to get the same accuracy.

\begin{example}\label{ex_2d}
{\rm We now provide a two-dimensional example to illustrate the accuracy of the proposed scheme:
\begin{flalign}\label{ex_2d_sys}
\left\{
\begin{aligned}
\left(\begin{aligned}
X_1(t)\\
X_2(t)
\end{aligned}\right)
=& \left(\begin{aligned}
X_1(0)\\
X_2(0)
\end{aligned}\right)
+ \int_0^t\left(\begin{aligned}
\alpha_1\sin^2(s+X_1(s))\\
\alpha_2\sin^2(s+X_2(s))
\end{aligned}\right)\di s\\
&+\int_0^t \left(\begin{aligned}
\alpha_2\cos^2(s+X_2(s))\\
\alpha_1\cos^2(s+X_1(s))
\end{aligned}\right) \di W_s,\\
\left(\begin{aligned}
Y_1(t)\\
Y_2(t)
\end{aligned}\right) =& \left(\begin{aligned}
\sin(T+X_1(T))\sin(T+X_2(T))\\
\cos(T+X_1(T))\cos(T+X_2(T))
\end{aligned}\right)\\
&+ \int_t^T \left(\begin{aligned}
f_s^1\\
f_s^2
\end{aligned}\right)\di s
-\int_t^T\left(\begin{aligned}
Z_1(s)\\
Z_2(s)
\end{aligned}\right) \di W_s,
\end{aligned}
\right.
\end{flalign}
where
$f^1_s, f^2_s$ are chosen as 
\begin{equation}
\begin{aligned}
f_s^1 =~& -(1+\alpha_1)\cos(s+X_1(s))\sin(t+X_2(s))-Z_2(s)\\
&+\frac12Y_1(s)(\alpha_2^2\cos^4(t+X_2(s))+\alpha_1^2\cos^4(s+X_1(s)))
-\alpha_1\alpha_2Y_2^3(s)\\
&-(1+\alpha_2)\sin(s+X_1(s))\cos(s+X_2(s)),\\
f_s^2 =~& (1+\alpha_1)\sin(s+X_1(s))\cos(s+X_2(s))-Z_1(s)\\
&+\frac12Y_2(s)(\alpha_2^2\cos^4(s+X_2(s))+\alpha_1^2\cos^4(s+X_1(s)))
-\alpha_1\alpha_2Y_1(s)Y_2^2(s)\\
&+(1+\alpha_2)\cos(s+X_1(s))\sin(s+X_2(s)),
\end{aligned}
\end{equation}
such that the exact solution of \eqref{ex_2d_sys} (with $\alpha_1 = \alpha_2 = 1/2$) is
\begin{equation}
\begin{aligned}
Y_1(t) =~&\sin(t+X_1(t))\sin(t+X_2(t)),\\
Y_2(t) =~&\cos(t+X_1(t))\cos(t+X_2(t)),\\
Z_1(t) =~& 0.5\cos(t+X_1(t))\sin(t+X_2(t))\cos^2(t+X_2(t))\\
&+0.5\sin(t+X_1(t))\cos(t+X_2(t))\cos^2(t+X_1(t)),\\
Z_2(t) =~& -0.5\sin(t+X_1(t))\cos^3(t+X_2(t))\\
&-0.5\cos^3(t+X_1(t))\sin(t+X_2(t)).
\end{aligned}
\end{equation}
The errors and convergence rates are listed in Table \ref{tab_ex_2d}, and
it is shown that the numerical Scheme
\ref{full_sch} is an order-$k$ scheme for
the two-dimensional example.

\begin{table}[H]\label{tab_ex_2d}
\caption{Errors and convergence rates for Example \ref{ex_2d}}
\begin{center}
\begin{tabular}{c|c|c|c|c|c|c|c}
\hline
Scheme                  &             &   N=8   &   N=16   &   N=32   &   N=64   &   N=128   & CR \\
\hline
\multirow{4}{*}{1-step} & $|Y_1(0)-Y_1^0|$ & 9.332E-02 & 4.701E-02 & 2.359E-02 & 1.182E-02 & 5.915E-03 & 0.995\\
                        \cline{2-8}
                        & $|Y_2(0)-Y_2^0|$ & 8.036E-02 & 4.110E-02 & 2.080E-02 & 1.047E-02 & 5.255E-03 & 0.984 \\
                        \cline{2-8}
                        & $|Z_1(0)-Z_1^0|$ & 6.061E-02 & 3.014E-02 & 1.496E-02 & 7.442E-03 & 3.711E-03 & 1.008 \\
                        \cline{2-8}
                        & $|Z_2(0)-Z_2^0|$ & 3.856E-02 & 1.763E-02 & 8.310E-03 & 4.027E-03 & 1.980E-03 & 1.069 \\
\hline
\multirow{4}{*}{2-step} & $|Y_1(0)-Y_1^0|$ & 4.235E-02 & 1.179E-02 & 3.093E-03 & 7.905E-04 & 1.998E-04 & 1.935 \\
                        \cline{2-8}
                        & $|Y_2(0)-Y_2^0|$ & 3.468E-02 & 9.672E-03 & 2.536E-03 & 6.480E-04 & 1.637E-04 & 1.935 \\
                        \cline{2-8}
                        & $|Z_1(0)-Z_1^0|$ & 3.680E-03 & 1.076E-03 & 2.892E-04 & 7.482E-05 & 1.902E-05 & 1.903 \\
                        \cline{2-8}
                        & $|Z_2(0)-Z_2^0|$ & 3.088E-03 & 6.779E-04 & 1.476E-04 & 3.428E-05 & 8.265E-06 & 2.139 \\
\hline
\multirow{4}{*}{3-step} & $|Y_1(0)-Y_1^0|$ & 1.244E-02 & 1.402E-03 & 1.621E-04 & 1.933E-05 & 2.354E-06 & 3.092 \\
                        \cline{2-8}
                        & $|Y_2(0)-Y_2^0|$ & 1.070E-02 & 1.211E-03 & 1.405E-04 & 1.680E-05 & 2.048E-06 & 3.087 \\
                        \cline{2-8}
                        & $|Z_1(0)-Z_1^0|$ & 3.802E-03 & 4.604E-04 & 5.627E-05 & 6.943E-06 & 8.618E-07 & 3.026 \\
                        \cline{2-8}
                        & $|Z_2(0)-Z_2^0|$ & 3.384E-03 & 3.573E-04 & 4.085E-05 & 4.849E-06 & 5.895E-07 & 3.117 \\
\hline
\multirow{4}{*}{4-step} & $|Y_1(0)-Y_1^0|$ & 1.525E-03 & 1.677E-04 & 1.267E-05 & 8.579E-07 & 5.561E-08 & 3.709 \\
                        \cline{2-8}
                        & $|Y_2(0)-Y_2^0|$ & 1.196E-03 & 1.356E-04 & 1.031E-05 & 6.996E-07 & 4.540E-08 & 3.696 \\
                        \cline{2-8}
                        & $|Z_1(0)-Z_1^0|$ & 5.720E-04 & 1.871E-05 & 6.615E-07 & 2.572E-08 & 1.133E-09 & 4.739 \\
                        \cline{2-8}
                        & $|Z_2(0)-Z_2^0|$ & 8.162E-05 & 1.231E-05 & 1.313E-06 & 9.819E-08 & 6.616E-09 & 3.415 \\
\hline
\end{tabular}
\end{center}
\end{table}

}

\end{example}


\subsection{Coupled cases}
In this subsection, we will test Scheme \ref{full_nsch} for solving
coupled FBSDEs. Compared to the numerical schemes for decoupled FBSDEs, an iterative process is
required for the coupled FBSDEs in Scheme \ref{full_nsch}. In our computations, it is noticed that this iterative process converges very quickly, and 3 or 4 iterations are enough up to the tolerance
$\epsilon_0 = 10^{-11}.$

\begin{example}\label{ex_c2} {\rm
Consider the following two systems of FBSDEs.
\begin{equation}\label{ex_c2_sys1}
\left\{\begin{aligned}
X_t = & x+\int_0^t \cos(s+X_s)(Y_s+Z_s) \di s
+ \int_0^t \sqrt{2}Y_s\sin(s+X_s)\di W_s,\\
Y_t =  & \sin(T+X_T) +\int_t^T (-\cos(s+X_s)-Y_s-Z_s\\
& \qquad\qquad \qquad \qquad +\sin^2(s+X_s)(Y_s+Z_s+Y_s^3))\di s
 - \int_t^T Z_s\di W_s.
\end{aligned}\right.
\end{equation}

\begin{equation}\label{ex_c2_sys2}
\left\{\begin{aligned}
X_t = & x + \int_0^t \cos(s+X_s)(Y_s+Z_s)\di s
+ \int_0^t \sqrt{2}(Y_s\sin(s+X_s)+1)\di W_s,\\
Y_t = & \sin(T+X_T)
+ \int_t^T (-\cos(s+X_s)-\cos^2(s+X_s)Z_s \\
&\qquad\qquad +(3\sin^2(s+X_s)+\sin^4(s+X_s))Y_s)\di s -\int_t^T Z_s\di W_s.
\end{aligned}
\right.
\end{equation}
The exact solutions $(Y_t, Z_t)$ of \eqref{ex_c2_sys1} and
\eqref{ex_c2_sys2} are respectively
\[
(\sin(t+X_t),  \sqrt{2}\cos(t+X_t)\sin^2(t+X_t))\text{ and }
(\sin(t+X_t),  \sqrt{2}\cos(t+X_t)(\sin^2(t+X_t)+1)).
\]
 Note that the diffusion coefficient
$\sigma_s= \sqrt{2}(Y_s\sin(s+X_s)+1)$ in \eqref{ex_c2_sys2}
satisfies condition
\eqref{sigma_ass}, but the diffusion coefficient
$\sigma= \sqrt{2}Y_s\sin(s+X_s)$ in \eqref{ex_c2_sys1} does not.
Moreover, these two diffusion coefficients
does not depend on $Z_s$.
The aim of this numerical example is to test
the ability of our Scheme \ref{full_nsch} to handle this kind of coupled FBSDE.
We choose $x=1.0$, and use Scheme \ref{full_nsch}
to solve problems \eqref{ex_c2_sys1} and \eqref{ex_c2_sys2}.
The errors, convergence rates and running
times are listed in Table \ref{nbsde_sys1} and Table \ref{nbsde_sys2}

\begin{table}[H]
\small
\begin{center}
\setlength{\belowcaptionskip}{0pt}
\caption{Errors and convergence rates for problem \eqref{ex_c2_sys1}.}\label{nbsde_sys1}
\begin{tabular}{c|c|c|c|c|c|c|c}
\hline
Scheme                  &             &   N=16   &   N=32   &   N=64   &   N=128  &   N=256  & CR         \\
\hline
\multirow{3}{*}{$k=1$} & $|Y^0-Y_0|$ & 4.648E-02 & 2.382E-02 & 1.205E-02 & 6.060E-03 & 3.039E-03 & 0.984 \\
                        \cline{2-8}
                        & $|Z^0-Z_0|$ & 9.148E-02 & 4.768E-02 & 2.421E-02 & 1.219E-02 & 6.123E-03 & 0.977 \\
                        \cline{2-8}
                        & RT          &   2.259   &   4.583   &   14.789  &   69.159  &  320.065  &       \\

\hline
\multirow{3}{*}{$k=2$} & $|Y^0-Y_0|$ & 1.914E-03 & 4.903E-04 & 1.242E-04 & 3.123E-05 & 7.828E-06 & 1.984 \\
                        \cline{2-8}
                        & $|Z^0-Z_0|$ & 3.690E-03 & 8.515E-04 & 2.023E-04 & 4.903E-05 & 1.204E-05 & 2.064 \\
                        \cline{2-8}
                        & RT          &   4.116   &   7.264   &   21.982  &   91.236  &  624.612  &       \\
\hline
\multirow{3}{*}{$k=3$} & $|Y^0-Y_0|$ & 1.168E-04 & 1.606E-05 & 2.143E-06 & 2.724E-07 & 3.433E-08 & 2.935 \\
                        \cline{2-8}
                        & $|Z^0-Z_0|$ & 1.066E-03 & 1.380E-04 & 1.752E-05 & 2.210E-06 & 2.794E-07 & 2.976 \\
                        \cline{2-8}
                        & RT          &   9.808   &   21.815  &   82.966  &  438.961  &  3762.623 &       \\
\hline
\multirow{3}{*}{$k=4$} & $|Y^0-Y_0|$ & 1.096E-05 & 6.493E-07 & 4.030E-08 & 2.545E-09 & 1.637E-10 &  4.006 \\
                        \cline{2-8}
                        & $|Z^0-Z_0|$ & 1.108E-04 & 6.481E-06 & 3.943E-07 & 2.380E-08 & 1.441E-09 &  4.055 \\
                        \cline{2-8}
                        & RT          & 13.020    & 38.194    & 179.564   & 1098.142  & 12302.685    &       \\
\hline
\end{tabular}

\vspace{0.3cm}
\setlength{\belowcaptionskip}{10pt}
\caption{Errors and convergence rates for problem \eqref{ex_c2_sys2}.}\label{nbsde_sys2}
\begin{tabular}{c|c|c|c|c|c|c|c}
\hline
Scheme                  &             &   N=16   &   N=32   &   N=64   &   N=128  &   N=256  & CR         \\
\hline
\multirow{3}{*}{$k=1$} & $|Y^0-Y_0|$ & 6.400E-02 & 2.857E-02 & 1.339E-02 & 6.465E-03 & 3.174E-03 & 1.081 \\
                        \cline{2-8}
                        & $|Z^0-Z_0|$ & 2.204E-01 & 1.078E-01 & 5.280E-02 & 2.605E-02 & 1.292E-02 & 1.023 \\
                        \cline{2-8}
                        & RT          & 2.712     & 4.577     & 12.907    & 62.917    & 395.805   &       \\
\hline
\multirow{3}{*}{$k=2$} & $|Y^0-Y_0|$ & 1.626E-03 & 2.933E-04 & 5.741E-05 & 1.257E-05 & 2.978E-06 & 2.273 \\
                        \cline{2-8}
                        & $|Z^0-Z_0|$ & 8.326E-02 & 2.260E-02 & 5.868E-03 & 1.496E-03 & 3.780E-04 & 1.948 \\
                        \cline{2-8}
                        & RT          & 3.791     & 6.697     & 26.820    & 137.501   & 1054.777   &       \\
\hline
\multirow{3}{*}{$k=3$} & $|Y^0-Y_0|$ & 1.697E-04 & 2.319E-05 & 3.027E-06 & 4.040E-07 & 5.245E-08 & 2.916 \\
                        \cline{2-8}
                        & $|Z^0-Z_0|$ & 1.642E-02 & 2.137E-03 & 2.740E-04 & 3.470E-05 & 4.361E-06 & 2.970 \\
                        \cline{2-8}
                        & RT          & 8.590     & 17.676    & 67.238    & 420.146   & 3362.073  &       \\
\hline
\multirow{3}{*}{$k=4$} & $|Y^0-Y_0|$ & 2.505E-05 & 2.573E-06 & 1.584E-07 & 1.024E-08 & 6.995E-10 & 3.822\\
                        \cline{2-8}
                        & $|Z^0-Z_0|$ & 1.131E-03 & 6.470E-05 & 3.285E-06 & 1.772E-07 & 9.845E-09 & 4.213 \\
                        \cline{2-8}
                        & RT          & 14.027    & 33.306    & 150.193   & 862.522   & 7046.404  &       \\
\hline
\end{tabular}
\end{center}
\end{table}
}
\end{example}

Tables \ref{nbsde_sys1} and \ref{nbsde_sys2} show that Scheme \ref{full_nsch}
is a high-order method (of order-$k$ accuracy) for solving coupled FBSDEs
with diffusion coefficient $\sigma$ independent of $Z_t$
and with the constrant \eqref{sigma_ass} or not,
and show that the scheme is
more efficient when the integer $k$ is bigger.

\begin{example}\label{ex_c4} {\rm To have a further investigation for Scheme \ref{full_nsch}, we now choose a diffusion coefficient
$\sigma_s$ that depends on $(X_s,Y_s,Z_s)$. The coupled FBSDEs are
\begin{equation}\label{ex_c4_eq}
\left\{
\begin{aligned}
X_t = & x -\int_0^t \frac12\sin(s+X_s)\cos(s+X_s)(Y_s^2+Z_s)\di s \\
& \quad+ \int_0^t \frac12\cos(s+X_s)(Y_s\sin(s+X_s)+Z_s+1)\di W_s,\\
Y_t = & \sin(T+X_T) + \int_t^T Y_sZ_s-\cos(s+X_s)\di s
- \int_t^T Z_s\di W_s.
\end{aligned}
\right.
\end{equation}
The exact solutions $Y_t$ and $Z_t$ of \eqref{ex_c4_eq} are
$Y_t= \sin(t+X_t)$ and $Z_t= \cos^2(t+X_t)$.

In this example, we choose $x=1.5.$ The errors, numerical convergence rates,
and running times are listed in Table \ref{tab_ex_c4}.

}

\begin{table}[H]\label{tab_ex_c4}
\small
\setlength{\belowcaptionskip}{0pt}
\caption{Errors and convergence rates for Example \ref{ex_c4}.}
\begin{center}
\begin{tabular}{c|c|c|c|c|c|c|c}
\hline
Scheme                  &             &   N=16   &   N=32   &   N=64   &   N=128  &   N=256  & CR         \\
\hline
\multirow{3}{*}{1-step} & $|Y^0-Y_0|$ & 1.747E-02 & 8.540E-03 & 4.222E-03 & 2.099E-03 & 1.047E-03 & 1.015 \\
                        \cline{2-8}
                        & $|Z^0-Z_0|$ & 2.067E-03 & 1.017E-03 & 5.042E-04 & 2.510E-04 & 1.252E-04 & 1.011 \\
                        \cline{2-8}
                        & RT          & 2.437     & 4.399     & 17.238    & 115.393   & 842.726   &       \\
\hline
\multirow{3}{*}{2-step} & $|Y^0-Y_0|$ & 3.146E-04 & 8.091E-05 & 2.055E-05 & 5.182E-06 & 1.299E-06 & 1.980\\
                        \cline{2-8}
                        & $|Z^0-Z_0|$ & 4.015E-05 & 9.504E-06 & 2.350E-06 & 5.781E-07 & 1.430E-07 & 2.031 \\
                        \cline{2-8}
                        & RT          & 8.018     & 14.680    & 56.705    & 263.602   & 2506.276  &       \\
\hline
\multirow{3}{*}{3-step} & $|Y^0-Y_0|$ & 4.467E-05 & 5.604E-06 & 6.973E-07 & 8.684E-08 & 1.083E-08 & 3.003 \\
                        \cline{2-8}
                        & $|Z^0-Z_0|$ & 5.351E-06 & 6.122E-07 & 7.233E-08 & 8.744E-09 & 1.075E-09 & 3.069 \\
                        \cline{2-8}
                        & RT          & 8.228     & 17.213    & 80.789    & 500.543   & 4205.121  &       \\
\hline
\multirow{3}{*}{4-step} & $|Y^0-Y_0|$ & 1.105E-06 & 6.914E-08 & 4.300E-09 & 2.683E-10 & 1.715E-11 & 3.996 \\
                        \cline{2-8}
                        & $|Z^0-Z_0|$ & 1.610E-07 & 1.123E-08 & 7.593E-10 & 4.901E-11 & 3.163E-12 & 3.911 \\
                        \cline{2-8}
                        & RT          & 13.982    & 40.198    & 182.862   & 1267.329  & 9419.386  &       \\
\hline
\end{tabular}
\end{center}
\end{table}

\end{example}

The results in Table \ref{tab_ex_c4} clearly show that Scheme \ref{full_nsch}
works well and is a order-$k$ numerical method
for solving the fully coupled FBSDEs in which the diffusion coefficient
 $\sigma$ depends on $X_t$, $Y_t$ and $Z_t$.

\section{Conclusions}\label{cons}

In this work, a new kind of
highly accurate multistep schemes
for solving fully coupled FBSDEs is proposed.
The key feature of such approaches is that the Euler method is used
to solve the forward SDE, which reduces dramatically the entire computational work,
and moreover, the numerical solution of the BSDE maintains the high-order accuracy.
Our numerical experiments show that the proposed multistep schemes
are effective and of high-order accuracy
for solving both decoupled and fully coupled FBSDEs.
We believe that the schemes proposed here are promising
in solving problems for example in finance, stochastic control, risk measure, etc. In our future studies, with proper updates for the spacial approaches, we will use our high-order schemes to solve higher dimensional FBSDEs.

\section*{Acknowledgments}
The authors would like to thank the associate editor and the referees for their valuable comments,
which improve much of the quality of the paper.

\end{document}